\documentclass[review]{elsarticle}
\usepackage[utf8]{inputenc}
\usepackage{amsmath,amsfonts}
\usepackage[top=1in,bottom=1in,left=1in,right=1in]{geometry}
\usepackage{xcolor}
\usepackage{multicol}
\usepackage{graphicx}
\usepackage[titlenumbered,linesnumbered,ruled,vlined]{algorithm2e}

\SetKwInput{KwInput}{Input}                % Set the Input
\SetKwInput{KwOutput}{Output}              % set the Output
\SetKwFor{Dor}{do}{}{endloop}

\graphicspath{
{figures/}
}

\usepackage{comment}
\usepackage[%
rm={oldstyle=false,proportional=true},%
sf={oldstyle=false,proportional=true},%
tt={oldstyle=false,proportional=true,variable=true},%
qt=false%
]{cfr-lm}
\usepackage{color}

\DeclareMathOperator{\diag}{diag}

\DeclareMathOperator{\trace}{trace}

\newcommand{\set}[1]{\left\{#1\right\}}
\newcommand{\xcap}[1]{X^{(#1)}}
\newcommand{\wcap}[1]{W^{(#1)}}
\newcommand{\pcap}[1]{P^{(#1)}}
\newcommand{\scap}[1]{S^{(#1)}}
\newcommand{\rcap}[1]{R^{(#1)}}
\newcommand{\ccap}[1]{C^{(#1)}}

\newcommand{\xcapt}[1]{X^{(#1)T}}

\newcommand{\scapt}[1]{S^{(#1)T}}

\newcommand{\thcap}[1]{\Theta^{(#1)}}

\newcommand{\nev}{n_{\mathrm{ev}}}

\newcommand{\lub}{\lambda_{\mathrm{ub}}}
\newcommand{\lf}{\lambda_{\mathrm{F}}}

\makeatletter
\def\verbatim@font{\fontfamily{cmtt}\scriptsize\selectfont}
\makeatother

\DeclareFontFamily{U}{MnSymbolC}{}
\DeclareSymbolFont{MnSyC}{U}{MnSymbolC}{m}{n}
\DeclareFontShape{U}{MnSymbolC}{m}{n}{
    <-6>  MnSymbolC5
   <6-7>  MnSymbolC6
   <7-8>  MnSymbolC7
   <8-9>  MnSymbolC8
   <9-10> MnSymbolC9
  <10-12> MnSymbolC10
  <12->   MnSymbolC12%
}{}
\DeclareMathSymbol{\powerset}{\mathord}{MnSyC}{180}

\usepackage[utf8]{inputenc}
\usepackage{fourier}
\usepackage{array}
\usepackage{booktabs}
\usepackage{tabularx}
\usepackage{makecell}
\usepackage{multirow}
\usepackage{tablefootnote}
\usepackage{comment}

\newcolumntype{Z}{ >{\centering\arraybackslash}X }

\title{Hybrid Eigensolvers for Nuclear Configuration Interaction Calculations}
\author[1]{Abdullah Alperen}
\author[1]{Metin Aktulga}
\author[2]{Pieter Maris}
\author[3]{Chao Yang}

\address[1]{Michigan State University, East Lansing, MI 48824, USA}
\address[2]{Dept. of Physics and Astronomy, Iowa State University, Ames, IA 50011, USA}
\address[3]{Lawrence Berkeley National Laboratory, CA 94720, USA}

\date{April, 2023}

\begin{document}

\begin{abstract}
We examine and compare several iterative methods for solving large-scale eigenvalue problems arising from nuclear structure calculations.  In particular, we discuss the possibility of using block Lanczos method, a Chebyshev filtering based subspace iterations and the residual minimization method accelerated by direct inversion of iterative subspace (RMM-DIIS) and describe how these algorithms compare with the standard Lanczos algorithm and the locally optimal block preconditioned conjugate gradient (LOBPCG) algorithm.  Although the RMM-DIIS method does not exhibit rapid convergence when the initial approximations to the desired eigenvectors are not sufficiently accurate, it can be effectively combined with either the block Lanczos or the LOBPCG method to yield a hybrid eigensolver that has several desirable properties.  We will describe a few practical issues that need to be addressed to make the hybrid solver efficient and robust. 
\end{abstract}

\maketitle

\section{Introduction}
The computational study of the structure of atomic nuclei involves solving the many-body Schr\"odinger equation for a nucleus consisting of $Z$ protons and $N$ neutrons, with $A=Z+N$ the total number of nucleons,
\begin{eqnarray}
    \hat{H} \; \Psi_i(\vec{r}_1, \ldots, \vec{r}_A) &=& E_i \, \Psi_i(\vec{r}_1, \ldots, \vec{r}_A) \,,
    \label{eq:Schrodinger}
\end{eqnarray}
where $\hat{H}$ is the nuclear Hamiltonian, $E_i$ are the discrete energy levels of the low-lying spectrum of the nucleus, and $\Psi_i$ the corresponding $A$-body wavefunctions.   A commonly used approach to address this problem is the no-core Configuration Interaction (CI) method (or No-Core Shell Model)~\cite{Barrett:2013nh}, in which the many-body Schr\"odinger equation, Eq.~(\ref{eq:Schrodinger}), becomes an eigenvalue problem
\begin{eqnarray}
    H \; x_i &=& \lambda_i \; x_i \,,
    \label{eq:eigen}
\end{eqnarray} 
where $H$ is an $n \times n$ square matrix that approximates the many-body Hamiltonian $\hat{H}$, $\lambda_i$ is the $i$th eigenvalue of $H$, and $x_i$ is the corresponding eigenvector.  The size $n$ of the symmetric matrix $H$ grows rapidly with the number of nucleons $A$ and with the desired numerical accuracy, and can easily be several billion or more; however, this matrix is extremely sparse, at least for nuclei with $A\ge 6$.  Furthermore, we are typically interested in only a few (5 to 10) eigenvalues at the low end of the spectrum of $H$.  An iterative method that can make use of an efficient Hamiltonian-vector multiplication procedure is therefore often the preferred method to solve Eq.~(\ref{eq:eigen}) for the lowest eigenpairs.

For a long time, the Lanczos algorithm~\cite{Lanczos1950} with full orthogonalization was the default algorithm to use because it is easy to implement and because it is quite robust even though it requires storing hundreds of Lanczos basis vectors.  Indeed, there are several software packages~\cite{Caurier1999,BROWN2014115,Johnson:2013bna,Johnson:2018hrx,Shimizu:2013xba,Shimizu:2019xcd} in which the Lanczos algorithm is implemented for nuclear structure calculations.  Here we focus on the software MFDn (Many-Fermion Dynamics for nuclear structure)~\cite{sternberg2008,MARIS201097,aktulga2014improving}, which is a hybrid MPI/OpenMPI code that is being used at several High-Performance Computing centers; it has recently also been ported to GPUs using OpenACC~\cite{WACCPD2021,MARIS2022101554}.  

In recent work~\cite{SHAO20181}, we have shown that the low-lying eigenvalues can be computed efficiently by using the Locally Optimal Block Preconditioned Conjugate Gradient (LOBPCG) algorithm~\cite{Knyazev2001}.  The advantages of the LOBPCG algorithm, which we will describe with some detail in the next section, over the Lanczos algorithm include
\begin{itemize}
\item The algorithm is a block method that allows us to multiply $H$ with several vectors simultaneously. That is, instead of an SpMV, one performs an Sparse Matrix-Matrix multiplication (SpMM) of a sparse square $n \times n$  matrix on a tall skinny $n \times n_b$ matrix at every iteration, which introduces an additional level of concurrency in the computation and enables us to exploit data locality better.  In order to converge 5 to 10 eigenpairs we typically use blocks of $n_b=8$ to $n_b=16$ vectors -- this can also be tuned to the hardware of the HPC platform. 
\item The algorithm allows us to make effective use of pre-existing approximations to several eigenvectors.
\item The algorithm allows us to take advantage of a preconditioner that can be used to accelerate convergence.
\item Other dense linear algebra operations can be implemented as level 3 BLAS.
%\item The algorithm generally has a lower memory requirement.
\end{itemize}
Even though Lanczos is efficient in terms of the number of Sparse Matrix Vector multiplications (SpMV) it uses, we have shown that the LOBPCG method often takes less wallclock time to run because performing a single SpMM on $n_b$ vectors is more efficient than performing $n_b$ SpMVs sequentially, which is required in the Lanczos algorithm.

However, there are occasionally some issues with LOBPCG:
\begin{itemize}
    \item The method can become unstable near convergence. Although methods for stabilizing the algorithm has been developed and implemented~\cite{hetmaniuk,robustlobpcg}, they do not completely eliminate the problem.
    \item Even though the algorithm in principle only requires storing three blocks of vectors, in practice, many more blocks of vectors are needed to avoid performing additional SpMMs in the Rayleigh-Ritz procedure. This is a problem for machines on which high bandwidth memory is in short supply (such as GPUs).
\end{itemize}

In this paper, we examine several alternative algorithms for solving large-scale eigenvalue problems in the context of nuclear configuration interaction calculations.  In particular, we will examine the block Lanczos algorithm~\cite{blockLanczos} and the Chebyshev Filtered Subspace Iteration (ChebFSI)~\cite{Saad_1984,Zhou_2014}.  Both are block algorithms that can benefit from an efficient implementation of the SpMM operation and can take advantage of good initial guesses to several eigenvectors, if they are available.  Neither one of these algorithms can incorporate a preconditioner, which is a main drawback. However, as we will show in Sect.~\ref{sec:examples}, in the early iterations of these algorithms, good approximations to the desired eigenpairs emerge quickly, even though the total number of SpMVs required to obtain accurate approximations can be higher compared to the Lanczos and LOBPCG algorithms.  This observation suggests that these algorithms can be combined with algorithms that are effective in refining existing eigenvector approximations. One such refinement algorithm is the Residual Minimization Method (RMM) with Direct Inversion of Iterative Subspace (DIIS) correction~\cite{zunger,jia,pulay80}. This algorithm has an additional feature that it can reach convergence to a specific eigenpair without performing orthogonalization against approximations to other eigenpairs as long as a sufficiently accurate initial guess is available.  Therefore, this algorithm can also be used to compute (or refine) different eigenpairs independently.  This feature introduces an additional level of concurrency in the eigenvalue computation that enhances the parallel scalability.

The paper is organized as follows. In the next section, we give an overview of the Lanczos, block Lanczos, LOBPCG, and ChebFSI algorithms.  We also describe the RMM-DIIS algorithm and discuss how it can be combined with (block) Lanczos, LOBPCG and ChebFSI to form a hybrid algorithm to efficiently compute the desired eigenpairs.  In Sect.~\ref{sec:examples}, we give several numerical examples to demonstrate the effectiveness of each of these algorithms in terms of the number of iterations.  The performance benefits of a hybrid algorithm designed by combining RMM-DIIS with one of the block algorithms are discussed in Sect.~\ref{sec:hybrid}. We also discuss the practical issue of deciding when switch to RMM-DIIS from block Lanczos or LOBPCG and how to implement RMM-DIIS to maximize its performance benefit.

\section{Numerical Algorithms}
\label{sec:algs}
We review several algorithms for computing a few algebraically smallest 
eigenvalues and the corresponding eigenvectors.
We denote the eigenvalues of the $n\times n$ nuclear CI Hamiltonian $H$ arranged
in an increasing order by $\lambda_1 \leq \lambda_2 \leq \cdots \leq \lambda_n$.
Their corresponding eigenvectors are denoted by $x_1$, $x_2$, ..., $x_n$.
We are interested in the first $\nev \ll n$ eigenvalues and eigenvectors. If we define
$X=[x_1,x_2,\dotsc,x_{\nev}]$ and
$\Lambda=\diag\set{\lambda_1,\lambda_2,\dotsc,\lambda_{\nev}}$, respectively, we have $HX=X\Lambda$.
For no-core CI calculations in nuclear physics we are often interested in only a few eigenvalues, that is $\nev \sim 5$ to $10$.  For the block solvers described below, it is generally beneficial to use slightly more vectors than the number of desired eigenvectors, $n_b > \nev$, and in practice we use $n_b = 8$ or $16$ for best performance.

\subsection{Lanczos Algorithm}
The Lanczos algorithm is a classical algorithm for solving large scale eigenvalue
problems. The algorithm generates an orthonormal basis of a $k$-dimensional 
Krylov subspace
\begin{equation}
\mathcal{K}(H;v_1) = \{ v_1, H v_1, ..., H^{k-1} v_1\},
\label{eq:krylov}
\end{equation}
where $v_1$ is an appropriately chosen and normalized starting guess. Such a basis is produced by a Gram--Schmidt process in which the key step of obtaining the $(j+1)$st basis vector $v_{j+1}$ is
\begin{equation}
   w_{j} = (I-V_{j} V_{j}^T) H v_{j}, \qquad v_{j+1} = w_{j} / \| w_{j} \|,
\label{eq:gs}
\end{equation}
where $V_{j}$ is a matrix that contains all previous orthonormal basis vectors, i.e.,
\[
V_{j} = \left( v_1, v_2, ..., v_{j}\right).
\]

The projection of $H$ into the $k$-dimensional subspace spanned by columns of $V_k$ is a tridiagonal matrix $T_k$ that satisfies
\begin{equation}
HV_k = V_k T_k + w_{k}e_k^T,
\label{eqn:lanfac}
\end{equation}
where $e_k$ is the last column of a $k\times k$ identity matrix. Approximate eigenpairs of $H$ are obtained by solving the $k\times k$ eigenvalue problem
\begin{equation}
T_k q = \theta q.
\label{eqn:trideig}
\end{equation}
It follows from \eqref{eqn:lanfac}, \eqref{eqn:trideig} and the fact that $V_k^T V_k=I_k$, $V_k^Tw_{k} = 0$ that the relative residual norm associated with an approximate eigenpair $(\theta, V_kq)$ can be estimated by
\begin{equation}
\frac{\| H (V_k q) - \theta (V_kq) \|}{|\theta|} 
= \frac{\|w_{k}\|\cdot | e_k^T q|}{|\theta|},
\label{eqn:ritzest1}
\end{equation}
for $\theta \neq 0$.

\subsection{Block Lanczos}

One of the drawbacks of the standard Lanczos algorithm is that it is not easy for this algorithm to take advantage of good initial guesses to more than one desired eigenvector.  Although we can take a simple linear combination (or average) of the initial guesses to several desired eigenvectors as the initial vector $v_1$, the Lanczos algorithm tends to converge to one of the eigenvectors much faster than others.

An algorithm that can take advantage of multiple starting guesses to different eigenvectors is the Block Lanczos algorithm.
The Block Lanczos algorithm generates an orthonormal basis of a block Krylov subspace 
\begin{equation}
\mathcal{K}(H;V_1) = \{ V_1, H V_1, ..., H^{k-1} V_1\},
\end{equation}
\label{eq:bkrylov}
where $V_1$ is a matrix that contains $n_b \ge \nev$ orthonormal basis vectors, where $\nev$ is the number of desired eigenvectors; in practice, we take $n_b$ slightly larger than $\nev$.  With this method we can make use of good initial guesses to the desired eigenvectors, e.g. obtained in smaller calculations.

The Gram--Schmidt process used to generate an orthonormal basis in Lanczos is replaced by a block Gram--Schmidt step that is characterized by
\begin{equation}
   W_{j} = (I-\mathbf{V}_{j} \mathbf{V}_{j}^T) H V_{j},
\label{eq:bgs}
\end{equation}
where the matrix $\mathbf{V}_{j}$ contains $j$ block orthonormal bases, i.e.,
\begin{equation}
\mathbf{V}_{j} = \left( V_1, V_2, ...,V_{j}\right).
\end{equation}
The normalization step in \eqref{eq:gs} is simply replaced by a QR factorization
step, i.e.
\[
W_{j} = V_{j+1}R_{j+1},
\]
where $V_{j+1}^T V_{j+1} = I_{n_b}$, and $R_{j+1}$ is an $n_b \times n_b$ upper triangular matrix.

The projection of $H$ into the subspace spanned by the columns of $\mathbf{V}_k$ is a block tridiagonal matrix $\mathbf{T}_k$ that satisfies
\begin{equation}
H\mathbf{V}_k = \mathbf{V}_k \mathbf{T}_k + W_{k}E_k^T,
\label{eqn:blanfac}
\end{equation}
where $E_k$ is the last $n_b$ columns of an $n_b\cdot k \times n_b\cdot k$ identity matrix.

Approximate eigenpairs of $H$ are obtained by solving the $n_b\cdot k\times n_b \cdot k$ eigenvalue problem
\begin{equation}
\mathbf{T}_k q = \theta q.
\label{eqn:btrideig}
\end{equation}
It follows from \eqref{eqn:blanfac}, \eqref{eqn:btrideig} and the fact that $\mathbf{V}_k^T \mathbf{V}_k=I_{n_bk}$, $\mathbf{V}_k^TW_{k} = \mathbf{0}$ that the relative residual norm associated with an approximate eigenpair $(\theta, V_kq)$ can be estimated by
\begin{equation}
\frac{\| H (V_k q) - \theta (V_kq) \|}{|\theta|} 
= \frac{\|W_{k}\|_F \cdot \| E_k^T q\|}{|\theta|},
\label{eqn:ritzest2}
\end{equation}
Algorithm~\ref{alg:blocklan} outlines the main steps of the Block Lanczos algorithm.

\begin{algorithm}[!tb]
\caption{The Block Lanczos algorithm}
\label{alg:blocklan}
\KwInput{The sparse matrix $H$, the number of desired eigenvalues $\nev$, an initial guess to the eigenvectors associated with the lowest $n_b \geq \nev$ eigenvalues $X^{(0)} \in \mathbb{R}^{n \times n_b}$,  convergence tolerance ($tol$) and maximum number of iteration allowed ($maxiter$);}
\KwOutput{$(\Lambda,X)$, where $\Lambda$ is a $\nev \times \nev$ diagonal matrix containing the desired eigenvalues, and $X \in \mathbb{R}^{n\times \nev}$ contains the corresponding eigenvector approximations;}

Generate $V_1 \in \mathbb{R}^{n\times n_b}$ that contains an orthonormal basis of $X^{(0)}$\;
$\mathbf{V}_1 = (V_1)$\;
$\mathbf{T}_1 = \mathbf{V}_1^T H \mathbf{V}_1$\;
Solve the projected eigenvalue problem $\mathbf{T}_{1} U = U \Theta$, where $U^TU = I$, $\Theta$ is a diagonal matrix containing eigenvalues of $\mathbf{T}_{1}$ in an ascending order\;
Determine number of converged eigenpairs $n_c$ by checking the Ritz residual estimate~\eqref{eqn:ritzest2}\;
goto 16 if $n_c \geq \nev$\;
\Dor{$i=1,2,\ldots,$maxiter}
{
    $W_{i} = (I-\mathbf{V}_{i} \mathbf{V}_{i}^T) H V_{i}$\;
    Generate $V_{i+1}$ that contains an orthonormal basis of  $W_{i}$\;
    $\mathbf{V}_{i+1} \leftarrow (\mathbf{V}_{i}\: V_{i})$\;
    Update $\mathbf{T}_{i+1} = \mathbf{V}_{i+1}^T H \mathbf{V}_{i+1}$\;
    Solve the projected eigenvalue problem $\mathbf{T}_{i+1} U = U \Theta$, where $U^TU = I$, $\Theta$ is a diagonal matrix containing eigenvalues of $\mathbf{T}_{i+1}$ in an ascending order\;
    $X_{i+1} = \mathbf{V}_{i+1} U(:,1:\nev)$\;
    Determine number of converged eigenpairs $n_c$ by checking the Ritz residual estimate~\eqref{eqn:ritzest2}\;
    $i\leftarrow i+1$ and exit the loop if $n_c \geq \nev$\;
}

$\Lambda \leftarrow \Theta$; $X \leftarrow X_i$\;
\end{algorithm}

Both the Lanczos and Block Lanczos algorithms produce approximation to the desired eigenvector in the form of 
\[
z = p_d(H) v_0,
\]
where $v_0$ is some starting vector and $d$ is the degree of the polynomial.  However, for the same number of multiplications of the sparse matrix $H$ with a vector (SpMVs), denoted by $m$, the degree of the polynomial generated in a Block Lanczos algorithm is $d = m/n_b$, whereas, in the standard Lanczos algorithm, the degree of the polynomial is $d = m$. Because the accuracy of the approximate eigenpairs obtained from the Lanczos and Block Lanczos methods is directly related to $d$, we expect more SpMVs to be used in a Block Lanczos algorithm to reach convergence. On the other hand, in a Block Lanczos method, one can perform $n_b$ SpMVs as a single SpMM on a tall skinny matrix consisting of a block of $n_b$ vectors, which is generally more efficient than performing $n_b$ separate SpMVs in succession. As a result, the Block Lanczos method can take less time even if it performs more SpMVs.

Furthermore, both with the standard Lanczos algorithm and with the Block Lanczos algorithm the memory to store the previous Lanczos vectors and the computational cost of the Gram--Schmidt process, Eq.~(\ref{eq:gs}), increase with the number of iterations.  In both the Lanczos and Block Lanczos algorithm, one can restart the algorithm after $d$ iteration with an improved set of approximate eigenvectors if the Gram--Schmidt process or the storage of $m$ Lanczos vectors becomes a bottleneck~\cite{arpack,krylovschur, ksblock,Shimizu:2019xcd}. This approach can also be used for check-point and restart purposes.

\subsection{LOBPCG}
It is well known that the invariant subspace associated with the smallest $\nev$ eigenvalues and spanned by columns of $X \in \mathbb{R}^{n\times \nev}$ is the solution to the trace minimization problem
\begin{equation}
\min_{X^TX = I} \trace(X^THX).
\label{eq:tracemin}
\end{equation}

The LOBPCG algorithm developed by Knyazev \cite{Knyazev2001} seeks to solve~\eqref{eq:tracemin}
by using the updating formula
\begin{equation}
\xcap{i+1} = \xcap{i}\ccap{i+1}_1+\wcap{i}\ccap{i+1}_2 + \pcap{i-1}\ccap{i+1}_3,
\label{eq:xupd}
\end{equation}
to approximate the eigenvector corresponding to the $\nev$ leftmost eigenvalues of $H$, where $\wcap{i} \in \mathbb{R}^{n\times \nev}$ is the preconditioned gradient of the Lagrangian
\begin{equation}
\mathcal{L}(X,\Lambda) = \frac12\trace(X^THX)
-\frac12\trace\left[(X^TX-I)\Lambda\right]
\label{eq:lag}
\end{equation}
associated with \eqref{eq:tracemin} at $\xcap{i}$, and $\pcap{i-1}$ is the search direction obtained in the $(i-1)$st iterate of the optimization procedure, and $\ccap{i+1}_1$, $\ccap{i+1}_2$, $\ccap{i+1}_3$ are a set of coefficient matrices of matching dimensions that are obtained by minimizing \eqref{eq:lag} within the subspace $\scap{i}$ spanned by
\begin{equation}
\scap{i} \equiv \left(
\xcap{i} \ \ \wcap{i} \ \ \pcap{i-1}
\right).
\label{eq:rrsubspace}
\end{equation}
To improve the convergence of the LOBPCG algorithm, one can include a few more vectors in $\xcap{i}$ so that the number of columns in $\xcap{i}$, $\wcap{i}$ and $\pcap{i}$ is $n_b \geq \nev$.

The preconditioned gradient $\wcap{i}$ can be computed as
\begin{equation}
\wcap{i} = K^{-1} (H\xcap{i} - \xcap{i}\thcap{i}),
\label{eq:wblock}
\end{equation}
where $\thcap{i} = \xcapt{i}H\xcap{i}$, and $K$ is a preconditioner that approximates $H$ in some way.  The subspace minimization problem that yields the coefficient matrix $\ccap{i+1}_1$, $\ccap{i+1}_2$, $\ccap{i+1}_3$, which are three block rows of a $3n_b\times n_b$ matrix $\ccap{i+1}$, can be solved as a generalized eigenvalue problem
\begin{equation}
\left( \scapt{i} H \scap{i} \right) \ccap{i+1}
= \left( \scapt{i}  \scap{i} \right) \ccap{i+1} D^{(i+1)},
\label{eq:rreig}
\end{equation}
where $D^{(i+1)}$ is a $n_b\times n_b$ diagonal matrix containing $n_b$ leftmost eigenvalues of the projected matrix pencil \\$\left( \scapt{i} H \scap{i}, \scapt{i}  \scap{i} \right)$.
The procedure that forms the projected matrices $\scapt{i} H \scap{i}$ and $\scapt{i}  \scap{i}$ and solves the projected eigenvalue problem~\eqref{eq:rreig} is often referred to
as the \emph{Rayleigh--Ritz} procedure~\cite{parlett}.
Note that the summation of the last two terms in \eqref{eq:xupd} represents the search direction followed in the $i$th iteration, i.e.,
\begin{equation}
\pcap{i} = \wcap{i}\ccap{i+1}_2 + \pcap{i-1} \ccap{i+1}_3.
\label{eq:pupd}
\end{equation}

\begin{algorithm}[!tb]
\caption{The basic LOBPCG algorithm}
\label{alg:lobpcg}
\KwInput{The sparse matrix $H$, a preconditioner $K$, an initial guess to the eigenvectors associated with the lowest $n_b \geq \nev$ eigenvalues $X^{(0)} \in \mathbb{R}^{n \times n_b}$, number of desired eigenvalues ($\nev$), convergence tolerance ($tol$) and maximum number of iteration allowed ($maxiter$);}
\KwOutput{$(\Lambda,X)$, where $\Lambda$ is a $\nev \times \nev$ diagonal matrix containing the desired eigenvalues, and $X \in \mathbb{R}^{n\times \nev}$ contains the corresponding eigenvector approximations;}

$[C^{(1)},\Theta^{(1)}]=\texttt{RayleighRitz}(H,\xcap{0})$\;
$\xcap{1} = \xcap{0} C^{(1)}$\;
$\rcap{1} = H X^{(1)} - X^{(1)} \Theta^{(1)}$\;
$P^{(0)} = \emptyset$\;
\Dor{$i=1,2,\ldots,$maxiter}
{
  $\wcap{i} = K^{-1} \rcap{i}$\;
  $\scap{i} = \left[ \xcap{i}, \wcap{i}, \pcap{i-1} \right]$\;
  $[\ccap{i+1}, \thcap{i+1}] = \texttt{RayleighRitz}(H,\scap{i})$\;
  $\xcap{i+1} = \scap{i} \ccap{i+1}$\;
  $\rcap{i+1} = H \xcap{i+1} - \xcap{i+1} \thcap{i+1}$\;
  $\pcap{i} = \wcap{i}\ccap{i+1}_2 + \pcap{i-1}\ccap{i+1}_3$\;
  Determine number of converged eigenpairs $n_c$ by comparing the relative norms of the leading $\nev$ columns of $\rcap{i+1}$ against the convergence tolerance $tol$\;
  exit if $n_c \geq \nev$\;
}
$\Lambda \leftarrow \Theta^{(i)}(1:\nev,1:\nev)$; $X \leftarrow X^{(i)}(:,1:\nev)$\;
\end{algorithm}

Algorithm~\ref{alg:lobpcg}  outlines the main steps of the basic LOBPCG algorithm.
The most computationally costly step of Algorithm~\ref{alg:lobpcg} is the multiplication of $H$ with a set of vectors.
Although it may appear that we need to perform such calculations in steps 8 (where the projected matrix $\scapt{i} H \scap{i}$ is formed) and 10, the multiplication of
$H$ with $\xcap{i}$, $\xcap{i+1}$ and $\pcap{i}$ can be avoided because $H\xcap{i+1}$ and $H\pcap{i}$ satisfy the following recurrence relationships
\begin{eqnarray}
 H\xcap{i+1} &=& H\xcap{i}\ccap{i+1}_1 + H\wcap{i}\ccap{i+1}_2 + H\pcap{i-1} \ccap{i+1}_3, \label{eq:upd1} \\
 H\pcap{i} &=& H\wcap{i}\ccap{i+1}_2 + H\pcap{i-1}\ccap{i+1}_3. \label{eq:upd2}
\end{eqnarray}
Therefore, the only SpMM we need to perform is $H\wcap{i}$.  
Again, for the nuclear CI calculations of interest, the dimension $n$ of the sparse symmetric matrix $H$ can be several billions, whereas $\wcap{i}$ is a tall skinny $n \times n_b$ matrix with $n_b$ typically of the order of 8 to 16.

\subsection{Chebyshev Filtering}
\begin{figure}[tb]
\centering
\includegraphics[width=0.7\textwidth]{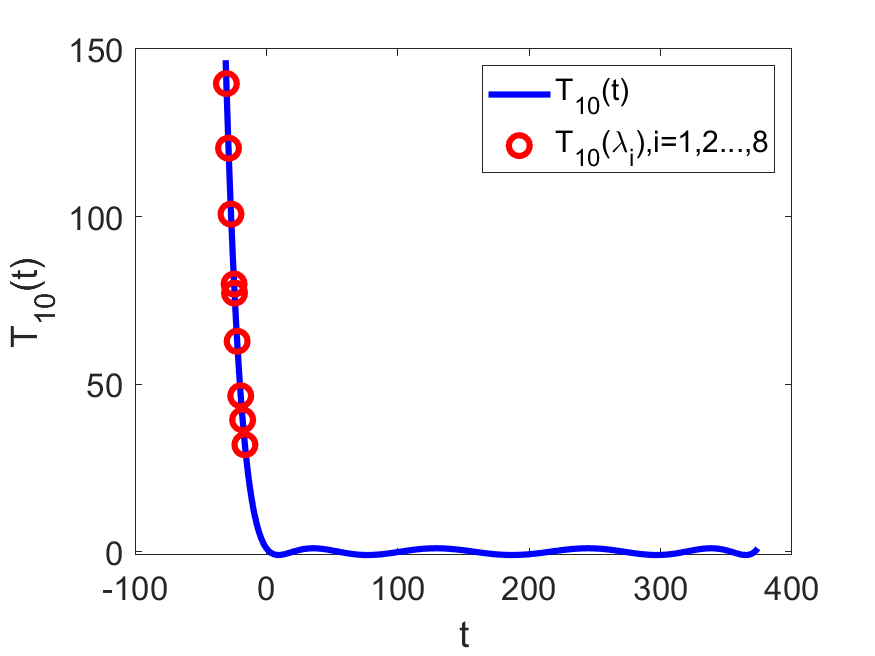}
\caption{Chebyshev polynomials of the first kind.}
\label{fig:Chebyshev_polynomial}
\end{figure}

An $m$th-degree Chebyshev polynomial of the first kind can be defined recursively as
\begin{equation}
T_m(t) = 2t T_{m-1}(t) - T_{m-2}(t),
\label{eq:chebypoly}
\end{equation}
with $T_0(t) = 1$ and $T_1(t) = t$.  The magnitude of $T_m(t)$ is bounded by 1 within $[-1,1]$ and grows rapidly outside of this interval. 
By mapping the unwanted eigenvalues of the nuclear many-body Hamiltonian $H$ enclosed by $[\lf,\lub]$ to $[-1,1]$ through the linear transformation $(t - c)/e$, where $c = (\lf+\lub)/2$ and $e = (\lub - \lf)/2$, we can use $\hat{T}_m(H) = T_m((H - cI)/e) v$ to amplify the eigenvector components in $v$ that correspond to eigenvalues outside of $[\lf,\lub]$.  Figure \ref{fig:Chebyshev_polynomial} shows a 10th degree Chebyshev polynomial defined on the spectrum of a Hamiltonian matrix and how the leftmost eigenvalues $\lambda_i$, $i=1,2,...,8$ are mapped to $T_{10}(\lambda_i)$.

Applying $T_m((H - cI)/e)$ repeatedly to a block of vectors $V$ filters out the eigenvectors associated with eigenvalues in $[\lf,\lub]$.  The desired eigenpairs can be obtained through the standard Rayleigh--Ritz procedure~\cite{parlett}. 

To obtain an accurate approximation to the desired eigenpairs, a high degree Chebyshev polynomial may be needed. Instead of applying a high degree polynomial once to a block of vectors, which can be numerically unstable, we apply Chebyshev polynomial filtering within a subspace iteration to iteratively improve approximations to the desired eigenpairs. We will refer to this algorithm as a Chebyshev Filtering based Subspace Iteration (ChebFSI).  The basic steps of this algorithm are listed in Algorithm~\ref{alg:filtering}.

\begin{algorithm}[!tb]
\caption{The Chebyshev filtering based subspace iteration (ChebFSI)}
\label{alg:filtering}
\KwInput{The sparse matrix $H$, an initial guess to the eigenvectors associated with the lowest $n_b \geq \nev$ eigenvalues $X^{(0)} \in \mathbb{R}^{n \times n_b}$, number of desired eigenvalues ($\nev$), the degree of the Chebyshev polynomial $d$; the spectrum cutoff $\lambda_{\mathrm{F}}$ and upper bound $\lambda_{\mathrm{ub}}$; convergence tolerance ($tol$) and maximum number of subspace iteration allowed ($maxiter$);}
\KwOutput{$(\Lambda,X)$, where $\Lambda$ is a $\nev \times \nev$ diagonal matrix containing the desired eigenvalues, and $X \in \mathbb{R}^{n\times \nev}$ contains the corresponding eigenvector approximations;}

$e = (\lub - \lf)/2$\;
$c = (\lub + \lf)/2$\;
$[Q,R]=\texttt{CholeskyQR}(H,\xcap{0})$\;
\Dor{$i=1,2,\ldots,$maxiter}
{
  $W = 2(HQ-cQ)/e$\;
  \Dor{$j=2,\ldots,d$}
  {
     $Y = 2(HW-cW)/e - Q$\;
     $Q\leftarrow W$\;
     $W\leftarrow Y$\;
  }
  $[Q,R]=\texttt{CholeskyQR}(Y)$\;
  $T = Q^THQ$\;
  Solve the eigenvalue problem $TS = S\Theta$, where $\Theta$ is diagonal, and update $Q$ by $Q\leftarrow QS$\;
  Determine number of converged eigenpairs $n_c$ by comparing the relative norms of the leading $\nev$ columns of $R = HQ - QD$ against the convergence tolerance $tol$\;
  exit if $n_c \geq \nev$\;
}
$\Lambda \leftarrow \Theta(1:\nev,1:\nev)$; $X \leftarrow Q(:,1:\nev)$\;
\end{algorithm}

Owing to the three-term recurrence in \eqref{eq:chebypoly}, $W=\hat{T}_m(H)V$ can be computed recursively without forming $\hat{T}_m(H)$ explicitly in advance. Lines 5 to 9 of Algorithm~\ref{alg:filtering} illustrates how this step is carried out in detail. To maintain numerical stability, we orthonormalize vectors in $W$. The orthonormalization can be performed by a (modified) Gram--Schmidt process or by a Householder transformation based QR factorization~\cite{Golub_1996}.

In Algorithm~\ref{alg:filtering}, the required inputs are a filter degree, $d$, an estimated upper bound of the spectrum of $H$, $\lub$, and an estimated spectrum cutoff level, $\lf$. The estimation of the upper bound $\lub$ can be calculated by running a few Lanczos iterations~\cite{vandetal2000,ZhouSaadTiagoEtAl2006,LiZhou2011}, and $\lf$ can often be set at 0 or an estimation of the $n_{b}+1$st leftmost eigenvalue of $H$ obtained from the Lanczos algorithm.  In the subsequent subspace iterations, $\lf$ can be modified based on more accurate approximations to the desired eigenvalues.  See the work of Saad~\cite{Saad_1984} and Zhou {\it et al.}~\cite{Zhou_2014} for more details on Chebyshev filtering.

\subsection{RMM-DIIS}
The Residual Minimization Method (RMM)~\cite{zunger,jia} accelerated by Direct Inversion of Iterative Subspace (DIIS)~\cite{pulay80} was developed in the electronic structure calculation community to solve a linearized Kohn-Sham eigenvalue problem in each self-consistency field (SCF) iteration. Given a set of initial guesses to the desired eigenvectors, $\{x_j^{0}\}$, $j=1,2,...,\nev$, the method produces successively more accurate approximations by seeking an optimal linear combination of previous approximations to the $j$th eigenvector by minimizing the norm of the corresponding sum of residuals. To be specific, 
let $x_j^{(i)}$, $i = 0,1,...,\ell-1$ be approximations to the $j$th eigenvector of $H$ obtained in the previous $\ell-1$ steps of the RMM-DIIS algorithm, and $\theta_j^{(i)}$ be the corresponding eigenvalue approximations. In the $\ell$th iteration (for $\ell > 1$), we first seek an approximation in the form of
\begin{equation}
\tilde{x}_j = \sum_{i=\min\{0,\ell-s\}}^{\ell-1} \alpha_i x_j^{(i)},
\label{eq:xsum}
\end{equation}
where 
\begin{equation}
\sum_{i=\ell-s}^{\ell-1} \alpha_i = 1, 
\label{eq:alphasum}
\end{equation}
for some fixed $1\leq s \leq s_{\max}$. The coefficients $\alpha_i$'s are obtained by 
solving the following constrained least squares problem
\begin{equation}
\min \| \sum_{i=\min\{\ell-s\}}^{\ell-1} \alpha_i r_j^{(i)} \|^2,
\label{eq:rmm}
\end{equation}
where $r_j^{(i)}=Hx_j^{(i)}-\theta_j^{(i)} x_j^{(i)}$ is the residual associated with the approximate eigenpair $(\theta_j^{(i)},x_j^{(i)})$, subject to the same constraint defined by \eqref{eq:alphasum}.
The constrained minimization problem can be turned into an unconstrained minimization problem by substituting $\alpha_{\ell-1} = 1 - \sum_{i=\min\{\ell-s\}}^{\ell-2} \alpha_i$
into \eqref{eq:rmm}.

Once we solve \eqref{eq:rmm}, we compute the corresponding residual
\[
\tilde{r}_j = H \tilde{x}_j - \tilde{\theta}_j \tilde{x}_j,
\]
where $\tilde{\theta}_j = \langle \tilde{x}_j, H \tilde{x}_j \rangle/\langle \tilde{x}_j, \tilde{x}_j \rangle$. A new approximation to the desired eigenvector is obtained by projecting $H$ into the two-dimensional subspace $W_j$ spanned by $\tilde{x}_j$ and $\tilde{r}_j$, and solving the $2\times 2$ generalized eigenvalue problem
\begin{equation}
(W_j^T H W_j) g = \theta (W_j^T W_j) g.
\label{eq:rr}
\end{equation}
Such an approximation can be written as
\begin{equation}
x_j^{(\ell)} = W_j g,
\label{eq:xuwg}
\end{equation}
where $g$ is the eigenvector associated with the smaller eigenvalue of the matrix pencil $(W_j^THW_j, W_j^TW_j)$.

Although Rayleigh--Ritz procedures defined by \eqref{eq:rr} and \eqref{eq:xuwg} are often used to compute the lowest eigenvalue of $H$, the additional constraint specified by \eqref{eq:xsum} and \eqref{eq:alphasum} keeps $x_j^{(\ell)}$ close to the initial guess of the $j$th eigenvector. Therefore, if the initial guess is sufficiently close to the $j$th eigenvector, $x_j^{(\ell)}$ can converge to this eigenvector instead of the eigenvector associated with the smallest eigenvalue of $H$.  Algorithm~\ref{alg:rmmdiis} outlines the main steps of the RMM-DIIS algorithm.

\begin{algorithm}[!tb]
\caption{The RMM-DIIS algorithm}
\label{alg:rmmdiis}
\KwInput{The sparse matrix $H$, an initial guess to the $\nev$ desired eigenvectors $\{x_j^{(0)}\}$, $j=1,2,...\nev$; maximum dimension of DIIS subspace $s$; convergence tolerance ($tol$) and maximum number of iteration allowed ($maxiter$);}
\KwOutput{$\{(\theta_j,x_j)\}$, $j=1,2,...,\nev$, where $\theta_j$ is the approximation to the $j$th lowest eigenvalue, and $x_j$ is the corresponding approximate eigenvector ;}

\For{j=1,2,\ldots,$\nev$}
{
   $x_j^{(0)} \leftarrow  x_j^{(0)}/\|x_j^{(0)}\|$\;
   $\theta^{(0)} = \langle x_j^{(0)}, Hx_j^{(0)}\rangle$\;
   $r_j^{(0)} = Hx_j^{(0)}-\theta_j^{(0)}x_j^{(0)}$\;
   $\tilde{x}_j^{(1)} = x_j^{(0)}$\;
   $\tilde{r}_j^{(1)} = r_j^{(0)}$\;
   \Dor{$i=1,2,\ldots,$maxiter}
   {
      \If{$i>1$}
      {
         Solve the residual minimization least squares problem \eqref{eq:rmm}\;
         Set $\tilde{x}_j^{(i)}$ according to \eqref{eq:xsum}\;
         $\tilde{x}_j^{(i)}\leftarrow \tilde{x}_j^{(i)}/\|\tilde{x}_j^{(i)}\|$\;
         Set the residual
         $\tilde{r}_j^{(i)} =  \sum_{\ell=\min\{0,i-s\}}^{i-1} \alpha_i x_j^{(\ell)}$\; 
      }
      
      Set $W_j = (\tilde{x}_j^{(i)}, \: \tilde{r}_j^{(i)})$\;
      Solve the Rayleigh--Ritz problem \eqref{eq:rr} and obtain $(\theta_j^{(i)},x_j^{(i)})$\;
      Compute the residual
         $r_j^{(i)} = Hx_j^{(i)} - \theta_j^{(i)}x_j^{(i)}$\;
      exit the do loop if $\|r_j^{(i)}\|/|\theta_j^{(i)}| < tol$\;
   }
}
\end{algorithm}

\subsection{Comparison summary}
The computational cost of all iterative methods discussed above is dominated by the the number of Hamiltonian matrix vector multiplications, that is, the number of SpMVs.  In the Lanczos algorithm, the number of SpMVs is the same as the number of iterations.  In block algorithms such as the Block Lanczos algorithm and the LOBPCG algorithm, the number of SpMVs is the product of the block size, $n_b$, and the number of iterations. The number of SpMVs used in the ChebFSI method is the product of the number of subspace iterations, the block size $n_b$, and the degree $d$ of the Chebyshev polynomial used.  The number of SpMVs used in RMM-DIIS is the sum of the RMM-DIIS iterations for each of the $\nev$ desired eigenpairs; depending on the architecture, one could combine $\nev$ SpMVs on single vectors in one SpMM on a tall skinny $n \times \nev$ matrix to improve performance.

In addition to SpMVs, some dense linear algebra operations are performed in these algorithms to orthonormalize basis vectors and to perform the Rayleigh--Ritz calculations. The cost of orthonormalization can become large if too many Lanczos iterations or Block Lanczos iterations are performed; for the Lanczos and Block Lanczos algorithm once can perform a restart once the orthonormalization cost becomes too large.
The orthonormalization cost is relatively small in LOBPCG, ChebFSI, and RMM-DIIS.

In Table~\ref{tab:memcomp}, we compare the memory usage of each method discussed above. Note that the first term for Lanczos, Block Lanczos, LOBPCG and ChebFSI in this table is generally the dominant term.  We use $\mathcal{O}(c)$ to denote a small multiple (i.e., typically 2 or 3) of $c$.
The number of iterations taken by a Block Lanczos iteration $k_{\mathrm{blockLan}}$ is typically smaller than the number of Lanczos iterations $k_{\mathrm{Lan}}$ when the same number of eigenpairs are computed by these methods. However, $k_{\mathrm{blockLan}}\cdot n_b$ is often larger than $k_{\mathrm{Lan}}$.
It is possible to use a smaller amount of memory in LOBPCG and ChebFSI at the cost of performing more SpMMs.  For example, if we were to explicitly compute $H\xcap{i+1}$ and $H\pcap{i}$ in the LOBPCG algorithm to perform the Rayleigh--Ritz calculation instead of updating these blocks according to \eqref{eq:upd1} and \eqref{eq:upd2} respective, we can reduce the LOBPCG memory usage to $4 n \cdot n_b + \mathcal{O}(9n_b^2)$.  For the RMM-DIIS algorithm, we assume that we compute one eigenpair at a time. The parameter $s_{\mathrm{max}}$ is the maximum dimension of the DIIS subspace constructed to correct an approximate eigenvector.  This parameter is often chosen to be between 10 and 20.
If we batch the refinement of several eigenvectors together to make use of SpMMs, the memory cost of RMM-DIIS will increase by a factor of $\nev$.

\begin{table}[htbp]
\small
\centering
\caption {A comparison of memory footprint associated with the Lanczos, Block Lanczos, LOBPCG, ChebFSI and RMM-DIIS methods.
} 
\label{tab:memcomp} 
\begin{tabular}{|c||c|}
\hline
\textbf{Method} & \textbf{Memory cost} \\ \hline\hline\\[-9pt]
Lanczos &  $n \cdot (k_{\mathrm{Lan}}+\nev) + \mathcal{O}\big(k_{\mathrm{Lan}}^2\big)$ \\ \hline\\[-9pt]
Block Lanczos & $n\cdot n_b\cdot (k_{\mathrm{block Lan}}+\nev) + \mathcal{O}\Big(\big(n_b\cdot k_{\mathrm{block Lan}}\big)^2\Big)$ \\ \hline\\[-9pt]
LOBPCG & $7n\cdot n_b + \mathcal{O}\big(9n_b^2\big)$ \\ \hline\\[-9pt]
ChebFSI & $4n\cdot n_b + \mathcal{O}\big(n_b^2\big)$ \\ \hline\\[-9pt]
RMM-DIIS & $n \cdot \big(3 \nev + s_\mathrm{max}\big)$ \\ \hline
\end{tabular}
\end{table}

\subsection{Hybrid Algorithms}

Among the methods discussed above, the Lanczos, Block Lanczos, LOBPCG, and ChebFSI methods can all proceed with an arbitrary starting guess of the desired eigenvectors although all, except the Lanczos algorithm, can (and generally do) benefit from the availability of good starting guesses to several eigenvectors.  On the other hand, as an eigenvector refinement method, the RMM-DIIS method requires a reasonably accurate approximation of the desired eigenvectors as a starting point.  Therefore, a more effective way to use the RMM-DIIS method is to combine it with one of the other methods, i.e., we can start with Lanczos, Block Lanczos, LOBPCG or ChebFSI method and switch to RMM-DIIS when the approximate eigenvectors become sufficiently accurate.  

In particular, in the Lanczos and Block Lanczos methods the basis orthogonalization cost as well as the memory requirement become progressively higher with increasing number of iterations.  Therefore, a notable benefit to switch from the Lanczos or Block Lanczos methods to RMM-DIIS is to lower the orthogonalization cost and memory requirement.

Although the orthogonalization cost and memory requirement for the LOBPCG method is fixed throughout all LOBPCG iterations, the subspace \eqref{eq:rrsubspace} from which eigenvalue and eigenvector approximations are drawn becomes progressively more ill-conditioned as the norms of the vectors in \eqref{eq:wblock} become smaller.  The ill-conditioned subspace can make the LOBPCG algorithm numerically unstable even after techniques proposed in~\cite{hetmaniuk,robustlobpcg} are applied.  Therefore, it may be desirable to switch from LOBPCG to RMM-DIIS, when the condition number of the subspace is not too large.  

The orthogonalization cost and memory requirement for ChebFSI are also fixed.  The method is generally more efficient in the early subspace iterations when $T_n(H)$ is applied to a block of vectors in each iteration.  However, as the approximate eigenvectors converge, applying $T_n(H)$ to a block of vectors in a single iteration results in a higher cost compared to RMM-DIIS that can refine each approximate eigenvector separately.  Again, it may become advantageous to switch to RMM-DIIS, when approximate eigenvectors become sufficiently accurate in ChebFSI.
\section{Numerical examples}
\label{sec:examples}
In this section, we compare and analyze the performance of each of the five algorithms presented in Section~\ref{sec:algs} using numerical examples, in terms of the number of SpMVs that are needed to achieve the requested tolerance for the lowest $\nev$ eigenpairs.  The initial tests were done using MATLAB, and we tested the effectiveness of algorithms in terms of the number of SpMVs required to reach the requested tolerance for $\nev$ eigenpairs.  Subsequently, we have also implemented these algorithms in Fortran90 and experimented with these algorithms in MFDn.

\subsection{Test problems}

The test problems we use are the many-body Hamiltonian matrices associated with four different nuclei, $^6$Li, $^7$Li, $^{11}$B, and $^{12}$C, where the superscripts indicate the total number of nucleons (protons plus neutrons) in the nuclei.  These Hamiltonian matrices are constructed in different CI model spaces labeled by the $N_{\max}$ parameter, using the two-body potential Daejeon16~\cite{Shirokov:2016ead}.  In Table~\ref{tab:mat_table}, we list the matrix size $n$ as well as the number of nonzero matrix elements in half the symmetric matrices.  Note that the matix size $n$ depends on $A$ (the number of nucleons) and the basis truncation parameter $N_{\max}$, but is independent of the interaction.  The number of nonzero matrix elements for a given nucleus and interaction is the same for any two-body interaction; but with three-body interactions, the number of nonzero matrix elements is more than an order of magnitude larger for the same matrix size.  Also, the number of iterations needed in any of the iterative solvers will generally depend on the interaction (as well as the nucleus and the truncation).

\begin{table}[htbp]
\centering
\caption{Test problems used in the numerical experiments.}
\label{tab:mat_table} 
\begin{tabular}{|c||c|r|r|}
\hline
\textbf{Nucleus} &  ${N_{\max}}$ & \textbf{Matrix size $n$} & \textbf{\# Non-zeros} \\ \hline\hline\\[-9pt]
$^6$Li   & 6 &   197,822 &  106,738,802 \\ \hline\\[-9pt]
$^7$Li   & 6 &   663,527 &  421,938,629 \\ \hline\\[-9pt]
$^{11}$B & 4 &   814,092 &  389,033,682 \\ \hline\\[-9pt]
$^{12}$C & 4 & 1,118,926 &  555,151,572 \\ \hline
\end{tabular}
\end{table}

Before solving eigenvalue problems for the nuclei and model spaces listed in Table~\ref{tab:mat_table}, we first construct a good initial guess for each of the $\nev$ desired eigenvectors by computing the lowest few eigenvalues and the corresponding eigenvectors of smaller Hamiltonian matrices constructed from a lower dimensional CI model space labelled by $N_{\max}-2$ values.  Table~\ref{tab:smaller_mat_table} shows the matrix dimensions $n$ and the number of nonzero matrix elements in these smaller Hamiltonians.  The initial guesses to the desired eigenvectors of the Hamiltonian matrices listed in Table~\ref{tab:mat_table} are obtained by padding the eigenvectors of the smaller Hamiltonian matrices by zeros to match the dimension of the original problems to be solved. As we can see, since the dimension of the problems listed in Table~\ref{tab:smaller_mat_table} are an order of magnitude or two smaller than the corresponding problems listed in Table~\ref{tab:mat_table}, they can be solved relatively easily and quickly by almost any method.

\begin{table}[htbp]
\small
\centering
\caption{The dimensions and number of nonzero matrix elements in half the matrices for $^6$Li, $^7$Li, $^{11}$B and $^{12}$C that are constructed in a lower dimensional configuration model space labelled by a smaller $N_{\max}$ value.}
\label{tab:smaller_mat_table} 
\begin{tabular}{|c||c|r|r|}
\hline
\textbf{System} & \textbf{$N_{\max}$} & \textbf{Matrix size $n$} & \textbf{\#Non-zeros} \\ \hline\hline\\[-9pt]
$^6$Li   & 4 &  17,040 &  4,122,448 \\ \hline\\[-9pt]
$^7$Li   & 4 &  48,917 & 14,664,723 \\ \hline\\[-9pt]
$^{11}$B & 2 &  16,097 &  2,977,735 \\ \hline\\[-9pt]
$^{12}$C & 2 &  17,725 &  3,365,099\\ \hline
\end{tabular}
\end{table}

All algorithms presented in section~\ref{sec:algs} have been implemented in MATLAB which is ideal for prototyping new algorithms. 

%We perform the numerical experiments presented below on a single AMD EPYC 7763 socket on a Perlmutter CPU node maintained at the National Energy Research Scientific Computing (NERSC) Center. The EPYC socket contains 64 cores and we disabled hyperthreading so 64 OpenMP threads were used.

For each test problem, we typically perform two sets of experiments for each algorithm. In the first set of experiments, we compute $\nev=5$ lowest eigenvalues and their corresponding eigenvectors. In the second set, we increase the number of eigenpairs to be computed to $\nev = 10$. All calculations are performed in double precision arithmetic.

\subsection{The performance of single method solvers}
In this section, we report and compare the performance of the Lanczos, block Lanczos, LOBPCG, ChebFSI and RMM-DIIS methods when they are applied to the test problems listed in Table~\ref{tab:mat_table}.  For block methods such as the block Lanczos, LOBPCG and ChebFSI methods, we set the block size, i.e., the number of vectors in the matrix $\mathbf{V}_j$ in \eqref{eq:bgs}, the matrix $X^{(i)}$ in \eqref{eq:xupd}, to $n_b = 8$ when computing the $\nev = 5$ lowest eigenpairs of $H$, or to $n_b = 16$ when computing the $\nev = 10$ lowest eigenpairs of $H$, and an SpMM is performed to multiply $H$ with $n_b$ vectors all at once, rather than performing $\nev$ separate SpMVs.  Even though RMM-DIIS is not a block method, the $\nev$ SpMVs performed in this algorithm can also be fused together as a single SpMM as we explain below.

For block methods, we choose the starting guess for each method as the eigenvectors of the Hamiltonian constructed in a smaller CI space (with a smaller $N_{\max}$ value), as listed in Table~\ref{tab:smaller_mat_table}, padded with zeros to match the size of the Hamiltonian in the larger CI space (with a larger $N_{\max}$ value) as mentioned earlier.  This is also used in the RMM-DIIS method which only requires a starting guess for each of the desired eigenpairs.  We should note that when such a starting guess is not sufficiently close to the desired eigenvector associated with the larger $N_{\max}$ value, the convergence of the RMM-DIIS method can be slow as we will see from the numerical examples presented below.  Note that it is also possible that the RMM-DIIS algorithm converges to a set of eigenpairs that do not correspond to the lowest $n_{\mathrm{ev}}$ eigenvalues of $H$.

For the Lanczos algorithm, we take the initial guess $v_0$ to be the linear combination of augmented eigenvectors associated with the lowest $\nev$ eigenvalues of the Hamiltonian constructed from the smaller CI space, i.e.,
\[
v_0 = \frac{1}{\nev} \sum_{i=1}^{\nev} \hat{z}_i \,,
\]
where $\hat{z}_i$ is the zero padded eigenvector associated with the $i$th eigenvalue of the Hamiltonian constructed from the smaller CI space. 

All methods are terminated when the relative residual norms or estimated residual norm associated with all desired eigenpairs are below the threshold of $\tau = 10^{-6}$. A relative residual norm for an approximate eigenpair $(\theta, z)$ is defined to be
\[
\frac{| H z - z \theta |}{|\theta|} \,.
\]
For the Lanczos and block Lanczos methods, we use \eqref{eqn:ritzest1} and \eqref{eqn:ritzest2} to estimate the relative residual norm without performing additional Hamiltonian matrix and vector multiplications.

The convergence of the ChebFSI algorithm depends on the choice of several parameters. Here, we use a 10th degree Chebyshev polynomial, i.e. $d=10$ in the ChebFSI method.  The upper bound of the spectrum $\lambda_{\mathrm{ub}}$ is determined by first running 10 Lanczos iterations and using Rayleigh--Ritz approximation to the largest eigenpairs ($\theta_{10}$, $u_{10}$) to set $\lambda_{\mathrm{ub}}$ to $\theta_{10} + \|r_{10}\|$, where $r_{10} = H u_{10} - \theta_{10}u_{10}$. We set the parameter $\lambda_{F}$ simply to 0 because the desired eigenvalues are bound states of the nucleus of interest and are expected to be negative. We apply the technique of deflation for converged eigenvectors, i.e., once the relative residual norm of an approximate eigenpair falls below the convergence tolerance of $10^{-6}$, we ``lock" the approximate eigenvector in place and do not apply $H$ to this vector in subsequent computations. These vectors will still participate in the Rayleigh--Ritz caculation performed in steps 11 and 12 of Algorithm~\ref{alg:filtering} and be updated as part of the Rayleigh--Ritz procedure.  

In Tables~\ref{tab:aggregated_5_eigs} and~\ref{tab:aggregated_10_eigs}, we compare the performance of Lanczos, block Lanczos, LOBPCG, ChebFSI and RMM-DIIS in terms of the total number SpMVs performed in each of these methods. It is clear from these tables that the Lanczos method uses the least number of SpMVs.
However, the number of SpMVs used by both the LOBPCG, block Lanczos and RMM-DIIS is within a factor of 3 when $\nev=5$ eigenpairs are computed.  Because 8 SpMVs can be fused as a single SpMM, which is more efficient, in the block Lanczos and LOBPCG method, the total wall clock time used by these methods can be less than that used by the Lanczos method.
Five SpMVs can also be fused in the RMM-DIIS method, even though the algorithm targets each eigenvalue separately. Because different eigenvalues may converge at a different rate, we may switch to using SpMVs when some of the eigenpairs converge.  
Whether it is beneficial to make such a switch depends on the performance difference between SpMM and SpMV, which may be architecture dependent. We will discuss this issue more in the next section. 

\begin{table}[htbp]
\small
\centering
\caption {SpMV count for different algorithms on MATLAB to compute five lowest eigenvalues.} 
\label{tab:aggregated_5_eigs} 
\begin{tabular}{|c||c|c|c|c|c|}
\hline
\textbf{System} & \textbf{Lanczos} & \textbf{Block Lanczos} & \textbf{LOBPCG} &  \textbf{ChebFSI} & \textbf{RMM-DIIS}\\ \hline
$^6$Li   &  95 & 208 & 184 & 480 & 174 \\ \hline
$^7$Li   & 109 & 280 & 240 & 960 & 291\\ \hline 
$^{11}$B &  82 & 240 & 192 & 950 & 152 \\ \hline
$^{12}$C & 106 & 248 & 192 & 890 & 181 \\ \hline
\end{tabular}
\end{table}

\begin{table}[!htbp]
\small
\centering
\caption {SpMV count for different algorithms on MATLAB to compute ten lowest eigenvalues.} 
\label{tab:aggregated_10_eigs} 
\begin{tabular}{|c||c|c|c|c|c|}
\hline
\textbf{System} & \textbf{Lanczos} & \textbf{Block Lanczos} & \textbf{LOBPCG} &  \textbf{ChebFSI} & \textbf{RMM-DIIS}\\ \hline
$^6$Li   & 114 & 464 & 464 &  690 & 686\\ \hline
$^7$Li   & 192 & 512 & 464 & 1,350 & 884 \\ \hline 
$^{11}$B & 180 & 480 & 432 & 2,470 &  $> 3,000$ \\ \hline 
%\tablefootnote{The relative residual norm associated with the 7th, 8th and 9th eigenvalues were still at $\approx 10^{-2}$ after 200 RMM-DIIS iterations.  \textcolor{red}{That is, after 2000 SpMVs?}} \\ \hline
$^{12}$C & 164 & 480 & 400 & 2,300 & 499\\ \hline
\end{tabular}
\end{table}

Table~\ref{tab:aggregated_10_eigs} shows that the number of SpMVs used in the Lanczos algorithm increases only slightly when we compute $\nev=10$ eigenpairs. However, the number of SpMVs required in the block Lanczos, LOBPCG, ChebFSI and RMM-DIIS increase at a higher rate.  This is mainly due to the fact that once a sufficiently large Krylov subspace has been constructed, we can easily obtain approximations to more eigenpairs without enlarging the subspace much further.  
Furthermore, RMM-DIIS fails to converge for the 7th, 8th, and 9th eigenpairs of $^{11}$B.  Inspection of the obtained spectra with the other methods reveals that these three eigenvalues for this particular case are near-degenerate.  Specifically, the 
eigenvalues of the 7th and 8th state differ by 0.6\% and these two states have the same conserved quantum numbers so they can easily mix; while the 8th and 9th states do have different quantum numbers, but their eigenvalues differ by only 0.03\%.  Indeed, it should not be surprising that refining individual eigenpairs may fail to converge for near-degenerate states.

\begin{figure}[!tb]
\centering
\includegraphics[scale=0.6]{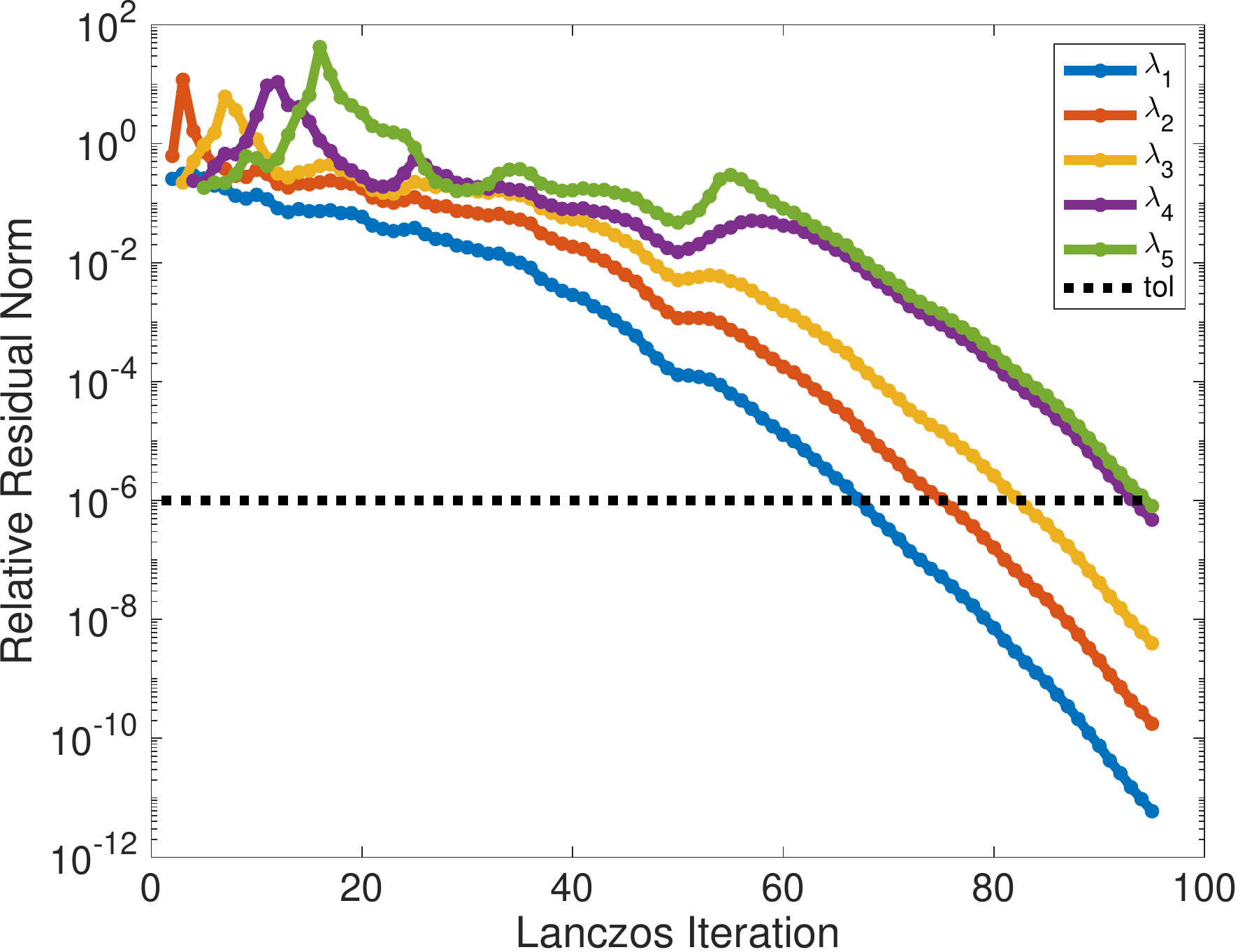}
\caption{The convergence of the lowest five eigenvalues of the $^6$Li Hamiltonian in the Lanczos algorithm as function of the number of iterations, which equals the number of SpMVs.}
\label{fig:Li6ConvLanczos}
\end{figure}
\begin{figure}[!tb]
\centering
\includegraphics[scale=0.6]{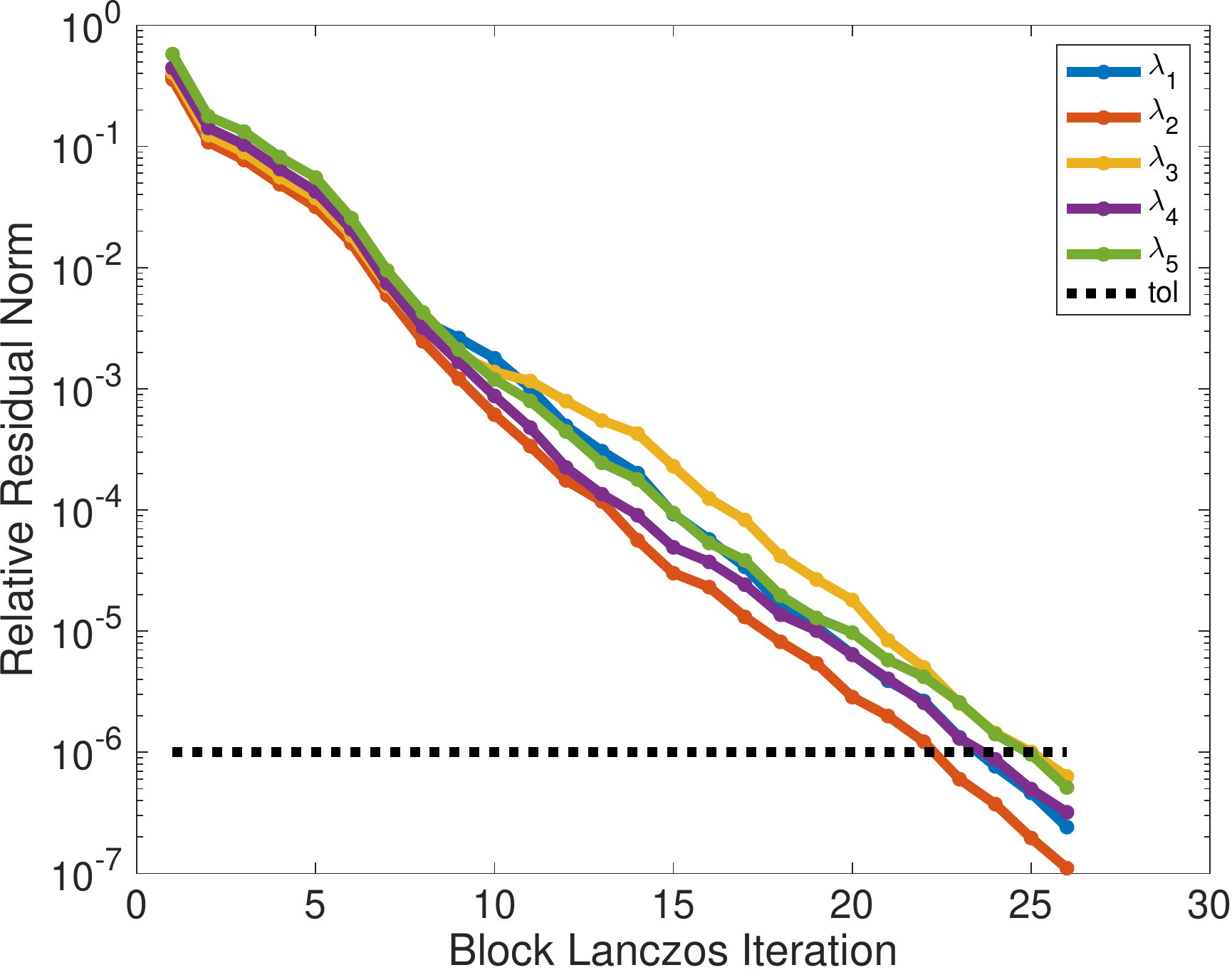}
\caption{The convergence of the lowest five eigenvalues of the $^6$Li Hamiltonian in the block Lanczos algorithm; 
the number of SpMVs is $n_b=8$ times the number of iterations.}
\label{fig:Li6ConvBlockLan}
\end{figure}
\begin{figure}[!tb]
\centering
\includegraphics[scale=0.6]{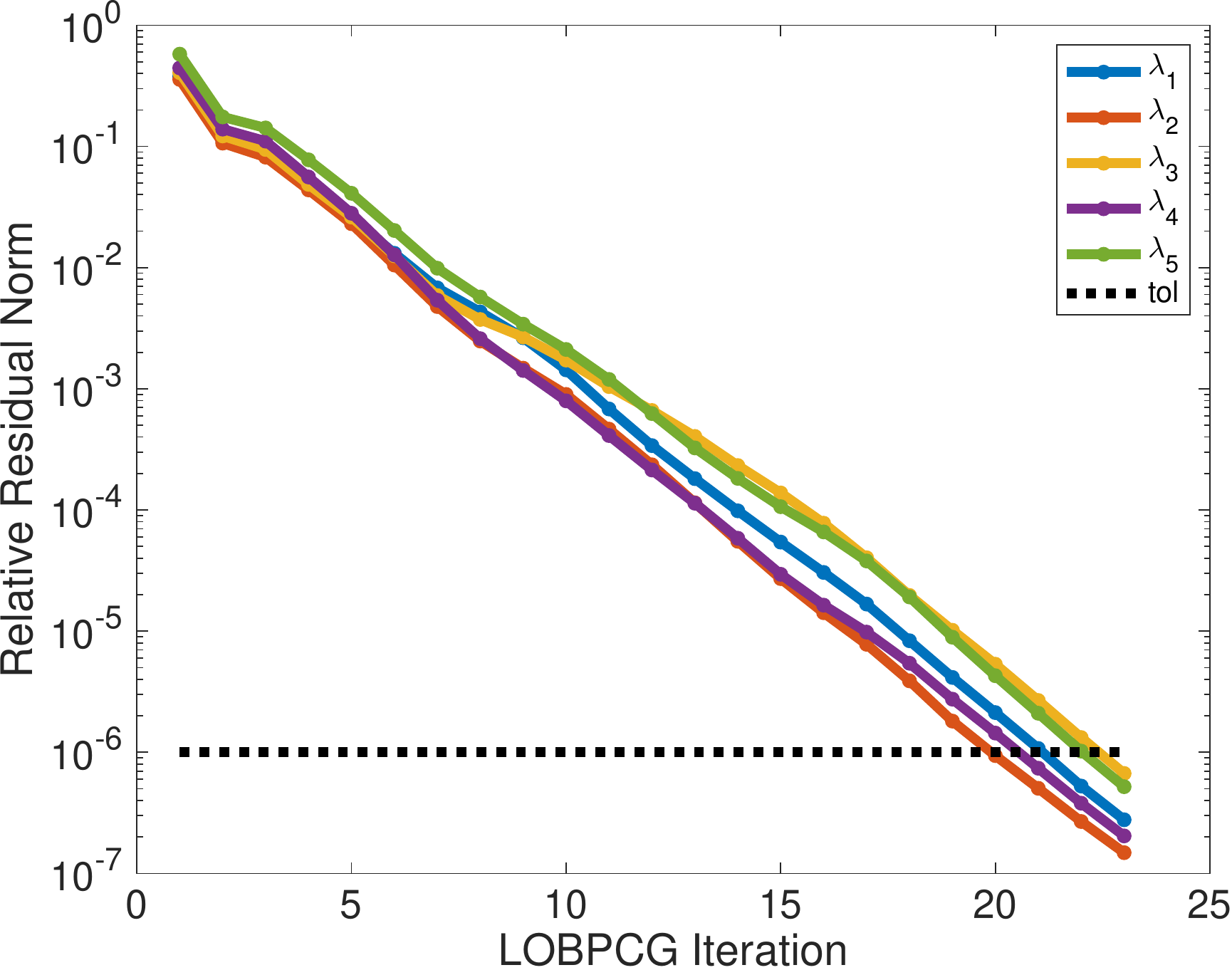}
\caption{The convergence of the lowest five eigenvalues of the $^6$Li Hamiltonian in the LOBPCG algorithm;
the number of SpMVs is $n_b=8$ times the number of iterations.}
\label{fig:Li6ConvLOBPCG}
\end{figure}
\begin{figure}[tb]
\centering
\includegraphics[scale=0.6]{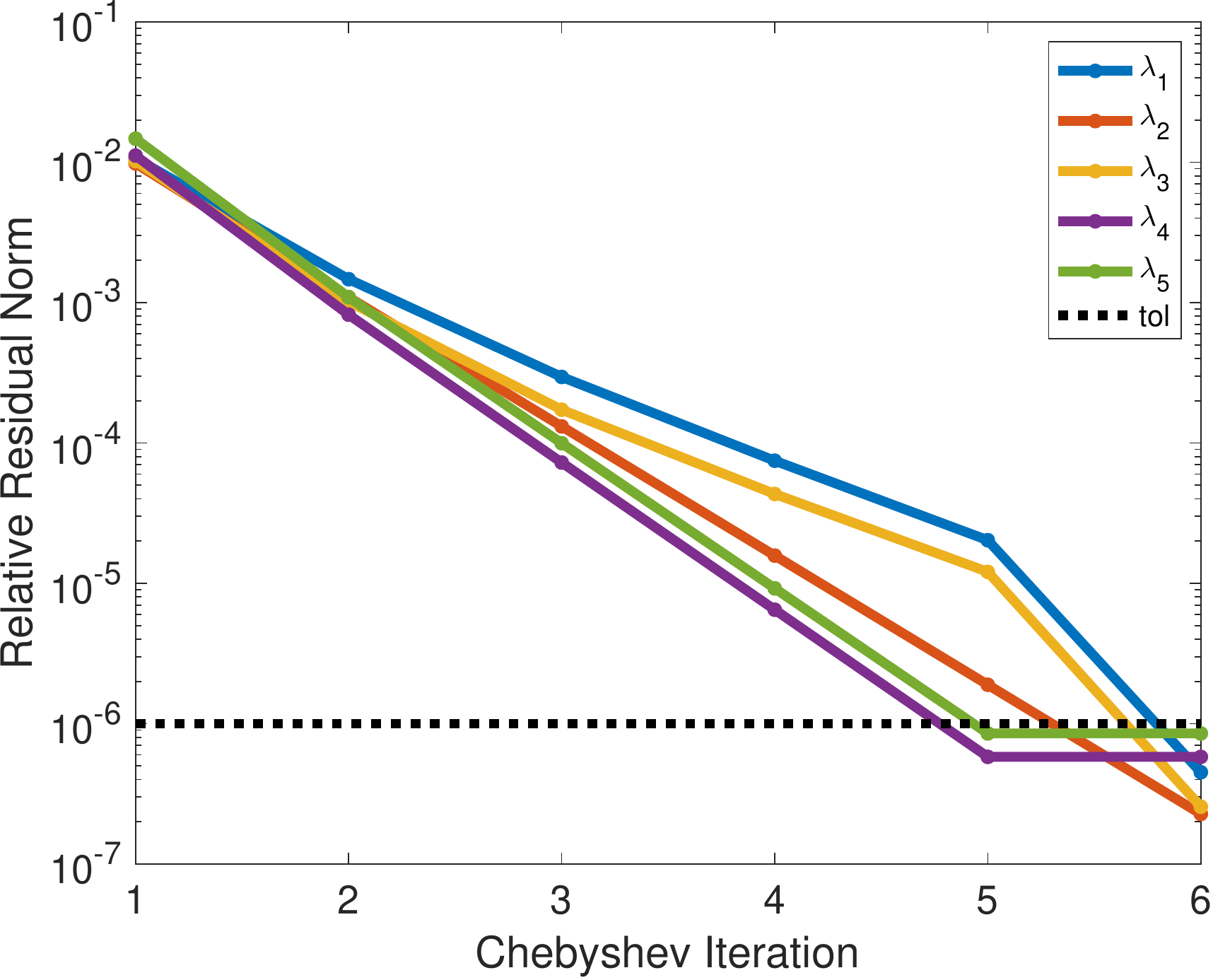}
\caption{The convergence of the lowest five eigenvalues of the $^6$Li Hamiltonian in the Chebyshev algorithm.
The number of SpMVs is $n_b=8$ times the degree of the Chebyshev polynomial, $d=10$, times the number of iterations.}
\label{fig:Li6ConvCheby}
\end{figure}
\begin{figure}[!tb]
\centering
\includegraphics[scale=0.6]{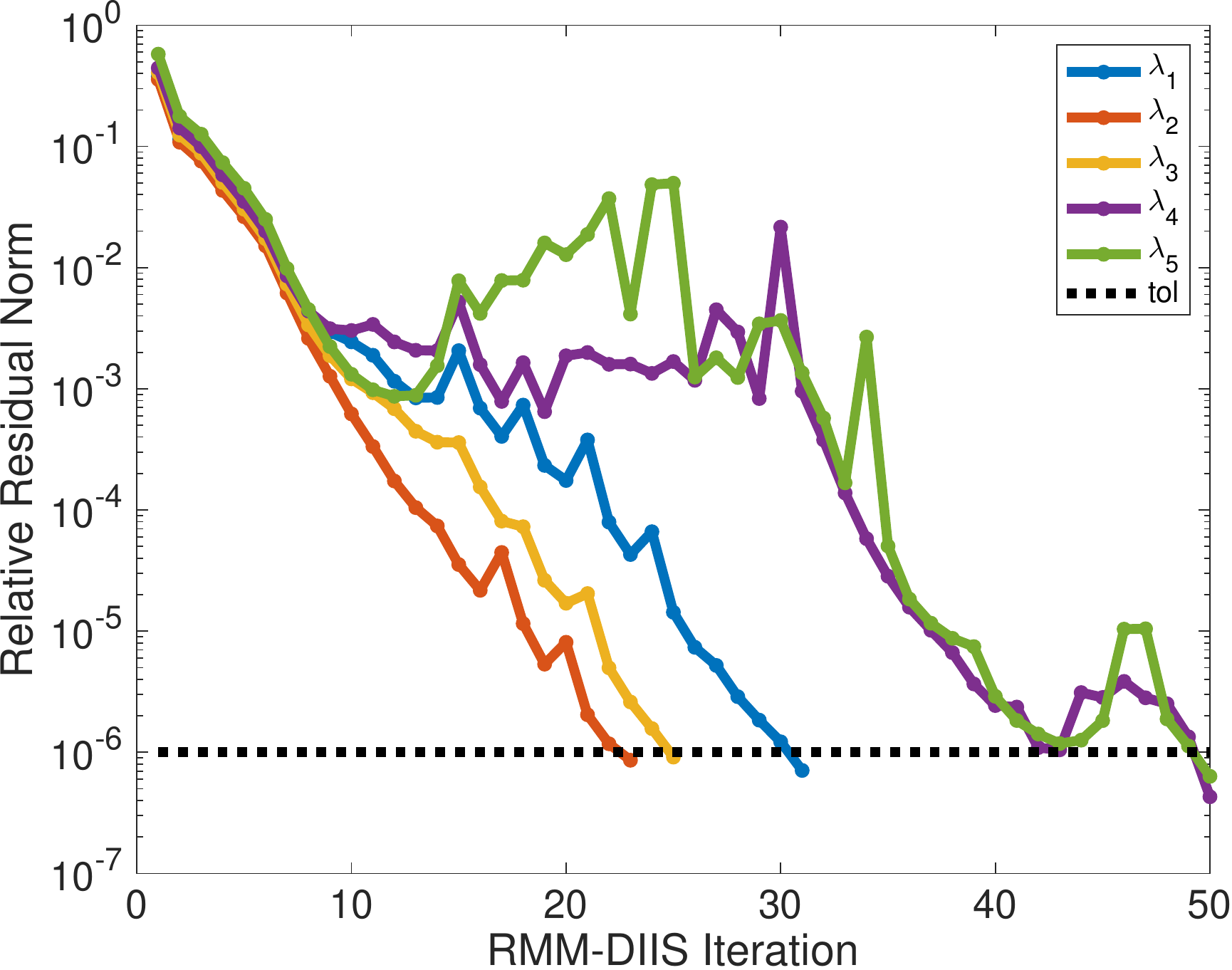}
\caption{The convergence of the lowest five eigenvalues of the $^6$Li Hamiltonian in the RMM-DIIS algorithm.
The number of SpMVs is the sum of the number of iterations for each eigenvalue.}
\label{fig:Li6ConvRMMDIIS}
\end{figure}

Figures~\ref{fig:Li6ConvLanczos},~\ref{fig:Li6ConvBlockLan},~\ref{fig:Li6ConvLOBPCG},~\ref{fig:Li6ConvCheby},~\ref{fig:Li6ConvRMMDIIS} show the convergence history of the Lanczos, block Lanczos, LOBCPG, ChebFSI and RMM-DIIS methods for $^6$Li.  In these figures, we plot the relative residual norm of each approximate eigenpair with respect to the iteration number.  We can clearly see that accurate approximations to some of the eigenpairs start to emerge in the Lanczos method when the dimension of the Krylov subspace (i.e., iteration number) is sufficiently large.  Note that the eigenpairs do not converge at the same rate.  Generally, the smallest eigenvalue converges first, followed by the second, third, fourth and the fifth eigenvalues.  However, these eigenvalues do not necessarily have to converge in order.  Although the relative residual for each eigenpair eventually goes below the convergence threshold of $10^{-6}$, the reduction of the relative residual norm is not monotonic with respect to the Lanczos iteration number.  The relative residual can sometimes increase after the Krylov subspace becomes sufficiently large and new spectral information becomes available.

All approximate eigenpairs appear to converge at a similar rate in the block Lanczos and LOBPCG methods.  This is one of the advantages of a block method.  The LOBPCG method performs slightly better in terms of the total number of SpMMs (or the number of iterations) used.  This is due to the fact that LOBPCG makes use of an effective preconditioner.

The ChebFSI is also a block method.  Figure~\ref{fig:Li6ConvCheby} shows that only five subspace iterations are required to obtain converged $\lambda_2$, $\lambda_4$ and $\lambda_5$, and two more subspace iterations are required to obtain converged $\lambda_1$ and $\lambda_3$.  The difference in the convergence rates for different eigenvalues is likely due to the variation in the contributions from different eigenvectors in the initial subspace constructed from the eigenvectors of the Hamiltonian associated with a smaller configuration space.  This could also be related to differences in the internal structure of the different states.  Although the number of subspace iterations used in ChebFSI is relatively small, each subspace iteration needs to perform $n_b \cdot d$ SpMVs, where $n_b$ is the number of vectors in the initial subspace and $d$ is the degree of the Chebyshev polynomial.
When $n_b=8$ and $d=10$, a total of 480 SpMVs are used to compute the lowest five eigenpairs as reported in Table~\ref{tab:aggregated_5_eigs}. Note that this count is less than $8\times 10 \times 7= 560$ because a deflation scheme that locks the converged eigenvector is used in ChebFSI.  When $n_b=16$, $d=10$, a total of 690 SpMVs are used to compute 10 lowest eigenpairs, as reported in Table~\ref{tab:aggregated_10_eigs}.  Because these SpMV counts are significantly higher than those used in other methods, ChebFSI appears to be not competitive for solving this type of eigenvalue problem. Therefore, from this point on, we will not discuss this method any further.

The convergence of the RMM-DIIS method is interesting. We observe from Figure~\ref{fig:Li6ConvRMMDIIS} that the first three eigenvalues of the $^6$Li Hamiltonian converge relatively quickly.  The number of RMM-DIIS iterations required to reach convergence is 31 for the first eigenvalue, 25 for the third eigenvalue and 23 for the second eigenvalue.  Altogether, 79 SpMVs are used to obtain accurate approximations to the three smallest eigenvalues and their corresponding eigenvectors, which is almost same as the 78 SpMVs used in the Lanczos algorithm for obtaining these three eigenpairs.  However, the fourth and fifth eigenvalues take much longer to converge. 
Yet, it is important to note that RMM-DIIS iterations that start with different initial guesses all converge to different eigenpairs even though no orthogonalization is performed between approximate eigenvectors produced in different RMM-DIIS iterations.  Table~\ref{tab:aggregated_10_eigs} shows that more RMM-DIIS iterations are needed to obtain accurate approximations to larger eigenvalues deeper inside the spectrum of $^6$Li Hamiltonian, partly due to the near-degeneracy of these eigenvalues in this case.
%\textcolor{red}{Remind me which interaction and which hbar-omega value you are using.  I suspect that the 4th and 5th eigenstates are both J=2 states, and are nearly degenerate.  They probably differ in their 'isospin', and depending on the interaction and hbar-omega value, I suspect that their may be a level crossing between the isospin-0 and isospin-2 state.  If that's indeed the case, we should mention that here.  Also indicate number of SpMVs in figure captions (to be done by PM, on Tuesday May 2) }
%\input discussion
\section{Performance of Hybrid Algorithms}
\label{sec:hybrid}
As we already discussed in Section~\ref{sec:algs}, the RMM-DIIS algorithm can be very effective when a good initial guess to the target eigenvector is available, with the caveat that it may perform poorly for near-degenerate states, even if they have different quantum numbers.  In principle, we can obtain reasonably accurate initial guesses for the target eigenvectors from any of the Lanczos, block Lanczos, LOBPCG or ChebFSI algorithms.  However, combining the RMM-DIIS algorithm with the Lanczos algorithm or the ChebFSI algorithm may be less attractive, partly because convergence rates for different eigenpairs in these algorithms are different, at least in the example discussed in Section~\ref{sec:examples}.  With the Lanczos algorithm we see from Figure~\ref{fig:Li6ConvLanczos} that the left most eigenvalues typically converge much faster than larger eigenvalues.  As a result, the RMM-DIIS method can only be effectively used to refine the eigenvectors associated with larger eigenvalues when they become sufficiently accurate.  At that point, the left most eigenvalues are likely to have converged already.  Similarly, in Figure~\ref{fig:Li6ConvCheby} we see that the first and third eigenpair converge slower than the other three eigenpairs.  Furthermore, the ChebFSI algorithm tends to be less efficient than either the block Lanczos or the LOBPCG algorithms; and also the Lanczos algorithm tends to be less efficient than block Lanczos or LOBPCG, at least when the multiplication of the sparse Hamiltonian $H$ with several vectors can be implemented efficiently.  We therefore do not consider the hybrid ChebFSI and RMM-DIIS nor the hybrid Lanczos and RMM-DIIS method here.

\begin{figure}[!tb]
\centering
\includegraphics[scale=0.6]{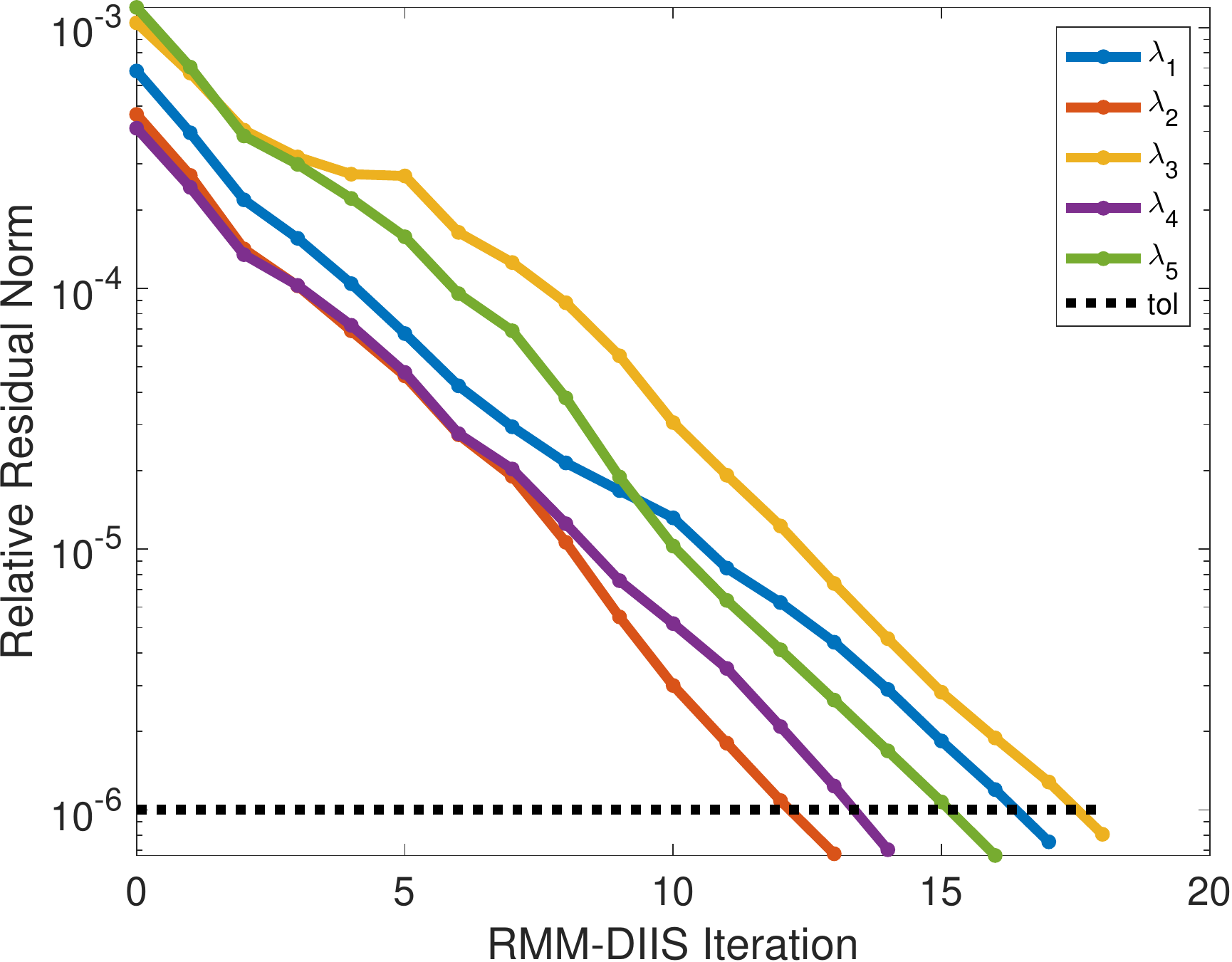}
\caption{The convergence of the lowest five eigenvalues of the $^6$Li Hamiltonian in the RMM-DIIS algorithm after 10 iterations of LOBPCG.}
\label{fig:Li6ConvLOBPCGRMMDIIS}
\end{figure}

In this section we therefore combine RMM-DIIS with either the block Lanczos or the LOBPCG algorithm to devise a hybrid algorithm that first runs a few block Lanczos or LOBPCG iterations and then use the RMM-DIIS method to refine approximated eigenvectors returned from the block Lanczos or the LOBPCG algorithm simultaneously.  Figure~\ref{fig:Li6ConvLOBPCGRMMDIIS} shows that RMM-DIIS indeed converges rapidly when the initial guesses to the eigenvectors associated with the lowest 5 eigenvalues are chosen to be the approximated eigenvectors produced from 10 LOBPCG iterations. 

\subsection{Switching to RMM-DIIS}
A practical question that arises when implementing a hybrid block Lanczos/RMM-DIIS or LOBPCG/RMM-DIIS eigensolver is how to decide when to switch from one algorithm to another.  Running more block Lanczos or LOBPCG iterations will produce more accurate approximations to the desired eigenvectors that can then be quickly refined by the RMM-DIIS algorithm.  But the cost of running block Lanczos or LOBPCG may dominate the overall cost of the computation.  On the other hand, running fewer block Lanczos or LOBPCG iterations may produce a set of approximate eigenvectors that require more RMM-DIIS iterations.  There is clearly a trade off, and the optimal point to switch from one to the other depends not only on the specific nucleus, model space, and interaction, but also on details of the numerical implementation and the computational hardware.

Table~\ref{tab:rmm-diis_lan_table_5_eigs} shows the total number of SpMVs required in a hybrid block Lanczos/RMM-DIIS algorithm in which 5, 10, 15 or 20 block Lanczos iterations are performed first, subsequently followed by the RMM-DIIS procedure.  We observe that this hybrid scheme uses fewer total SpMVs for all test problems than that used in a simple block Lanczos or RMM-DIIS run.  A similar observation can be made from Table~\ref{tab:rmm-diis_lan_table_10_eigs} when 10 lowest eigenvalues are computed except for one case, out of 16 cases: for $^7$Li, the hybrid algorithm that starts with 5 block Lanczos iterations uses more SpMVs than the non-hybrid block Lanczos only algorithm.  The optimal number of block Lanczos iterations we should perform (in terms of the total SpMV count) before switching to the RMM-DIIS procedure is problem dependent.  Because good approximations to interior eigenvalues only emerge when the Krylov subspace produced by the block Lanczos iteration is sufficiently large, more block Lanczos iterations are required in the hybrid algorithm to yield an optimal SpMV count when 10 eigenvalues and the corresponding eigenvectors are needed, as we can see in Table~\ref{tab:rmm-diis_lan_table_10_eigs}.
\begin{table}[!tb]
\small
\centering
\caption {SpMV count in hybrid block Lanczos/RMM-DIIS for computing five lowest eigenvalues.} 
\label{tab:rmm-diis_lan_table_5_eigs} 
\begin{tabular}{|c||c|c|c|c|c|}
\hline
{\textbf{Nucleus}} 
 & \textbf{5 block Lanczos} & \textbf{10 block Lanczos} & \textbf{15 block Lanczos} & \textbf{20 block Lanczos} & \textbf{block Lanczos only}\\ \hline
$^6$Li   & 183 & 163 & 173 & 185 & 208\\ \hline
$^7$Li   & 205 & 218 & 234 & 221 & 280\\ \hline 
$^{11}$B & 170 & 178 & 195 & 199 & 240 \\ \hline
$^{12}$C & 179 & 180 & 179 & 201 & 248 \\ \hline
\end{tabular}
\end{table}
\begin{table}[!tb]
\small
\centering
\caption {SpMV count in hybrid block Lanczos/RMM-DIIS for computing 10 lowest eigenvalues.} 
\label{tab:rmm-diis_lan_table_10_eigs} 
\begin{tabular}{|c||c|c|c|c|c|}
\hline
{\textbf{Nucleus}} 
 & \textbf{5 block Lanczos} & \textbf{10 block Lanczos} & \textbf{15 block Lanczos} & \textbf{20 block Lanczos} & \textbf{block Lanczos only}\\ \hline
$^6$Li   & 422 & 361 & 347 & 367 & 464 \\ \hline
$^7$Li   & 544 & 422 & 418 & 402 & 512 \\ \hline 
$^{11}$B & 419 & 422 & 392 & 382 & 480 \\ \hline
$^{12}$C & 430 & 415 & 369 & 369 & 480 \\ \hline
\end{tabular}
\end{table}

Also a hybrid LOBPCG and RMM-DIIS algorithm performs better than the LOBPCG algorithm by itself, provided a sufficiently large number of LOBPCG iterations are performed first to obtain good starting vectors for the subsequent RMM-DIIS procedure, as can be seen in Tables~\ref{tab:rmm-diis_lobpcg_table_5_eigs} and \ref{tab:rmm-diis_lobpcg_table_10_eigs}. 
\begin{table}[!tb]
\small
\centering
\caption {SpMV count in hybrid LOBPCG/RMM-DIIS for computing five lowest eigenvalues.} 
\label{tab:rmm-diis_lobpcg_table_5_eigs} 
\begin{tabular}{|c||c|c|c|c|c|}
\hline
{\textbf{System}} 
 & \textbf{5 LOBPCG} & \textbf{10 LOBPCG} & \textbf{15 LOBPCG} & \textbf{20 LOBPCG} & \textbf{LOBPCG only}\\ \hline
$^6$Li   & 191 & 158 & 163 & 165 & 184 \\ \hline
$^7$Li   & 211 & 200 & 190 & 194 & 240 \\ \hline 
$^{11}$B & 260 & 162 & 159 & 165 & 192 \\ \hline
$^{12}$C & 186 & 157 & 153 & 166 & 192 \\ \hline
\end{tabular}
\end{table}
\begin{table}[!tb]
\small
\centering
\caption {SpMV count in hybrid LOBPCG/RMM-DIIS for computing 10 lowest eigenvalues.} 
\label{tab:rmm-diis_lobpcg_table_10_eigs} 
\begin{tabular}{|c||c|c|c|c|c|}
\hline
{\textbf{System}} 
 & \textbf{5 LOBPCG} & \textbf{10 LOBPCG} & \textbf{15 LOBPCG} & \textbf{20 LOBPCG} & \textbf{LOBPCG only}\\ \hline
$^6$Li & 432 & 344 & 337 & 356 & 464 \\ \hline
$^7$Li & 529 & 398 & 364 & 365 & 464 \\ \hline 
$^{11}$B & 758 & 397 & 334 & 355 & 423 \\ \hline
$^{12}$C & 523 & 340 & 317 & 334 & 400 \\ \hline
\end{tabular}
\end{table}
From Tables~\ref{tab:rmm-diis_lan_table_5_eigs},~\ref{tab:rmm-diis_lan_table_10_eigs},~\ref{tab:rmm-diis_lobpcg_table_5_eigs} and~\ref{tab:rmm-diis_lobpcg_table_10_eigs}, we observe that the total SpMV count is optimal when the number of block Lanczos or LOBPCG iterations is sufficiently large, but not too large.  However, the optimal number of block Lanczos/LOBPCG iterations is problem dependent.

One practical way to determine when to switch from block Lanczos or LOBPCG to RMM-DIIS is to examine the change in the approximate eigenvalue.  Because RMM-DIIS can be viewed as an eigenvector refinement method, it works well when an approximate eigenvalue becomes sufficiently accurate while the corresponding approximate eigenvector still needs to be corrected.  Since we do not know the exact eigenvalues in advance, we use the average changes in the relative difference between approximate eigenvalues obtained in two consecutive iterations (e.g., the $(k-1)$st and the $k$th iterations) defined as
\begin{equation}
\tau = \frac{1}{n_{\mathrm{ev}}} \sqrt{  \sum_{j=1}^{n_{\mathrm{ev}}} {\left(\frac{\theta_j^{(k)} -\theta_j^{(k-1)}}{\theta_j^{(k)}}\right)^2}}.
\label{eq:switchtau}
\end{equation}
as a metric to decide when to switch from block Lancos or LOBPCG to RMM-DIIS.

Tables~\ref{tab:optimal_5_eigs} and~\ref{tab:optimal_10_eigs} show the optimal SpMV count of the hybrid block Lanczos/RMM-DIIS and LOBPCG/RMM-DIIS algorithms when computing five or $10$ lowest eigenvalues, respectively.  Unlike Tables~\ref{tab:rmm-diis_lan_table_5_eigs},~\ref{tab:rmm-diis_lan_table_10_eigs},~\ref{tab:rmm-diis_lobpcg_table_5_eigs} and~\ref{tab:rmm-diis_lobpcg_table_10_eigs} where we show the results only for $5, 10, 15$ or $20$ iterations, these optimal numbers are found by exhaustively trying every possible iteration count of block Lanczos and LOBPCG before switching to RMM-DIIS.  Along with the optimal numbers, Tables~\ref{tab:optimal_5_eigs} and~\ref{tab:optimal_10_eigs} also show the SpMV count of the hybrid algorithms when block Lanczos and LOBPCG procedures are stopped as $\tau$ goes below $10^{-4}$ and $10^{-7}$.

Table~\ref{tab:optimal_10_eigs} shows that the number of SpMV counts in both the hybrid block Lanczos/RMM-DIIS and LOBPCG/RMM-DIIS algorithms appear to be nearly optimal when $\tau \leq 10^{-7}$.  However, when fewer eigenvalues are needed, we can possibly stop the block Lanczos procedure in the hybrid block Lanczos/RMM-DIIS algorithm when $\tau$ is below $10^{-4}$ as shown in Table~\ref{tab:optimal_5_eigs}.  The $\tau \leq 10^{-7}$ criterion still seems to be a good one for the hybrid LOBPCG/RMM-DIIS algorithm even when fewer eigenpairs are needed.

\begin{table}[!tb]
\small
\centering
\caption {SpMV count for hybrid algorithms wrt switching threshold on MATLAB to compute five lowest eigenvalues.}
\begin{tabular}{|c||c|c|c|c|c|c|}
\hline
\textbf{Nucleus} & \textbf{Opt Hybrid Block Lan} & $\tau \leq 10^{-4}$ & $\tau \leq 10^{-7}$ & \textbf{Opt Hybrid LOBPCG} & $\tau \leq 10^{-4}$ & $\tau \leq 10^{-7}$ \\ \hline
$^6$Li   & 157 & 194 & 169 & 156 & 187 & 159 \\ \hline
$^7$Li   & 205 & 211 & 224 & 190 & 210 & 190 \\ \hline 
$^{11}$B & 166 & 176 & 195 & 158 & 188 & 159 \\ \hline
$^{12}$C & 172 & 180 & 179 & 152 & 169 & 153 \\ \hline
\end{tabular}
\label{tab:optimal_5_eigs} 
\end{table}
\begin{table}[!tb]
\small
\centering
\caption {SpMV count for hybrid algorithms wrt switching threshold on MATLAB to compute 10 lowest eigenvalues.}
\begin{tabular}{|c||c|c|c|c|c|c|}
\hline
\textbf{Nucleus} & \textbf{Opt Hybrid Block Lan} & $\tau \leq 10^{-4}$ & $\tau \leq 10^{-7}$ & \textbf{Opt Hybrid LOBPCG} & $\tau \leq 10^{-4}$ & $\tau \leq 10^{-7}$ \\ \hline
$^6$Li   & 340 & 384 & 346 & 330 & 393 & 335 \\ \hline
$^7$Li   & 401 & 515 & 418 & 363 & 478 & 363 \\ \hline 
$^{11}$B & 379 & 525 & 384 & 332 & 625 & 333 \\ \hline
$^{12}$C & 362 & 408 & 368 & 312 & 437 & 317 \\ \hline
\end{tabular}
\label{tab:optimal_10_eigs} 
\end{table}

\subsection{Optimal implementation}
As we indicated earlier, the SpMV count may not be the best metric to evaluate the performance of a block eigensolver that performs SpMVs for a block of vectors as a SpMM operation.  It has also been shown in \cite{aktulga2014optimizing,SHAO20181} that on many-core processors, the LOBPCG algorithm can outperform the Lanczos algorithm even though its SpMV count is much higher.  The reason that a block algorithm can outperform a single vector algorithm such as the Lanczos method on a many-core processor is that SpMM can take advantage of high concurrency and memory locality.  

The performance of SpMV and SpMM can be measured in terms of number of floating point operations performed per second (GFLOPS).  In the following experiments, 
we measure the performance these computational kernels by running the the LOBPCG solver implemented in the MFDn software on a single AMD EPYC 7763 socket on a Perlmutter CPU node maintained at the National Energy Research Scientific Computing (NERSC) Center. The EPYC socket contains 64 cores~\cite{nersc}. We disabled hyperthreading so that 64 OpenMP threads were used within a single MPI rank.
Tables~\ref{tab:gflops_epyc_rhs8} and ~\ref{tab:gflops_epyc_rhs16} show that the SpMM GFLOPS measured in the LOBPCG algorithm is much higher than the SpMV GFLOPS measured in the Lanczos algorithm.  As a result, we can evaluate the performance of a block eigensolver by dividing the actual SpMV count by the ratio between SpMM and SpMV GFLOPS to obtain an ``effective" SpMV count.  For example, because the SpMM and SpMV GFLOPS ratio is $24.2/4.1\approx 5.9$ for $^6$Li, the effective SpMV count for performing 23 iterations of the LOBPCG algorithm with a block size 8 is $23 \times 8 / 5.9 \approx 31$ which is much lower than the 95 SpMVs performed in the Lanczos method, even though the actual number of SpMVs performed in the LOBPCG iteration is $23\times 8 = 184 > 95$. 
%\textcolor{red}{This is all still in MATLAB? The performance in GFLOPS is likely to be different in MFDn, which is F90, using a non-standard SpMV/SpMM implementation, and for large runs (and we are interested in large runs), data movement becomes more of an issue -- a small run might fit entirely in cache, but that's certainly not the case for medium size runs; and in particular, in large runs one should also take the MPI overhead into consideration.  Typically I do runs for matrix dimensions that are between 1,000 and 10,000 times the dimension of the matrices considered here -- the FLOP comparisons made here are, unfortunately, completely irrelevant for those large-scale runs, and we may want to shorten this discussion a bit...}
%
\begin{table}[htbp]
\centering
\caption {GFLOPs achieved by the SpMV/SpMM kernels within MFDn on one AMD EPYC node with 64 threads in a single MPI rank. The sparse Hamiltonian is applied to 8 vectors in SpMM.} 
\label{tab:gflops_epyc_rhs8} 
\begin{tabular}{|c||c|c|c|}
\hline
\textbf{Nucleus} & \textbf{SpMV GFLOPS} & \textbf{SpMM GFLOPS} & \textbf{Ratio} \\ \hline
$^6$Li   & 4.1 & 24.2 & 5.9\\ \hline
$^7$Li   & 6.7 & 36.0 & 5.4 \\ \hline
$^{11}$B & 6.4 & 33.5 & 5.2\\ \hline
$^{12}$C & 6.0 & 28.0 & 4.7 \\ \hline
\end{tabular}
\end{table}

\begin{table}[htbp]
\centering
\caption {GFLOPs achieved by the SpMV/SpMM kernels within MFDn on one AMD EPYC node with 64 threads. The sparse Hamiltonian is applied to 16 vectors in SpMM.} 
\label{tab:gflops_epyc_rhs16} 
\begin{tabular}{|c||c|c|c|}
\hline
\textbf{Nucleus} & \textbf{SpMV GFLOPS} & \textbf{SpMM GFLOPS} & \textbf{Ratio} \\ \hline
$^6$Li   & 4.1 & 38.9 & 9.5 \\ \hline
$^7$Li   & 6.7 & 56.0 & 8.4 \\ \hline
$^{11}$B & 6.4 & 50.3 & 7.9 \\ \hline
$^{12}$C & 6.0 & 47.5 & 7.9 \\ \hline
\end{tabular}
\end{table}

Strictly speaking, RMM-DIIS is a single vector method, i.e., in each RMM-DIIS run, we refine one specific eigenvector associated with a target eigenvalue which has become sufficiently accurate as discussed earlier.  However, because the refinement of different eigenvectors can be performed independently from each other, we can perform several RMM-DIIS refinements simultaneously. The simultaneous RMM-DIIS runs can be implemented by batching the SpMVs in each RMM-DIIS iteration together as a single SpMM. This step constitutes the major cost of the RMM-DIIS method. The least squares problems given in Eqs.~(\ref{eq:rmm})--(\ref{eq:xuwg}) for different eigenvectors can be solved in sequence in each step, since they do not cost much computation.

Because different eigenvectors may converge at different rates as we have already seen in Figure~\ref{fig:Li6ConvRMMDIIS}, we need to decide what to do when one or a few eigenvectors have converged.  One possibility is to decouple the batched RMM-DIIS method after a certain number of eigenvectors have converged, and switch from using a single SpMM in a coupled RMM-DIIS implementation to using several SpMVs in a decoupled RMM-DIIS implementation. Another possibility is to just keep using the coupled RMM-DIIS with a single SpMM in each step without updating the eigenvectors that have already converged. In this case, the SpMM calculation performs more floating point operations than necessary. However, because an SpMM can be carried out at a much higher GFLOPs than an SpMV, the overall performance of the computation may not be degraded even with the extra computation.

\begin{figure}[htbp]
\centering
\includegraphics[scale=0.6]{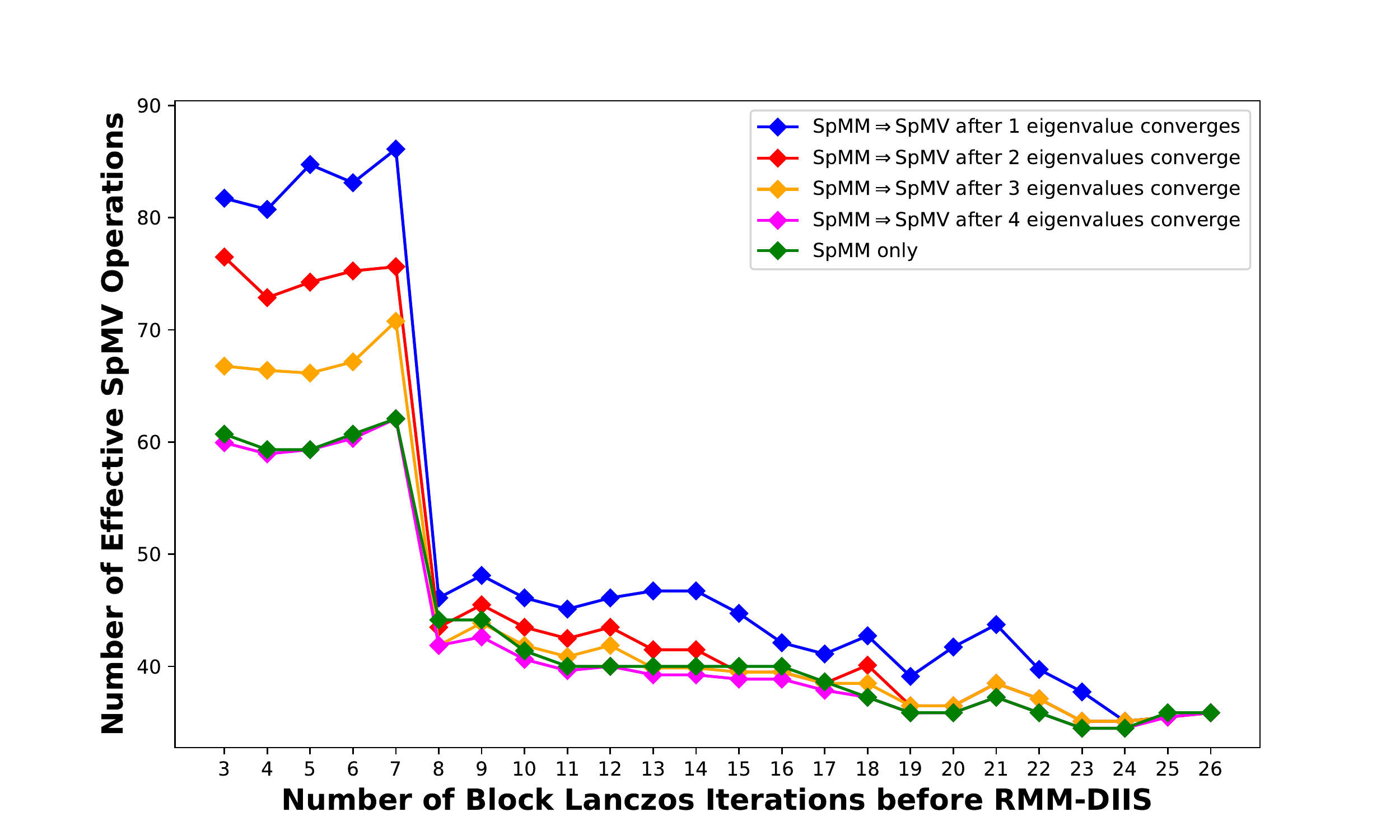}
\caption{Total effective SpMV cost of the block Lanczos/RMM-DIIS algorithm to compute the lowest five eigenvalues for $^6$Li w.r.t. the switching strategy from SpMM to SpMV within RMM-DIIS.}
\label{fig:switch_blocklan_li6}
\end{figure}
In Figure~\ref{fig:switch_blocklan_li6}, we show the effective number of SpMVs performed in several hybrid block Lanczos and RMM-DIIS runs for the $^6$Li test problem.  The horizontal axis represents the number of block Lanczos iterations performed before we switch to RMM-DIIS. The blue dots show the number of effective SpMVs used in the hybrid method when we switch from SpMM to SpMV after one eigenvector has converged. Similarly, the red, orange and magenta dots show the effective SpMV counts when we switch from SpMM to SpMV after two, three and four eigenvectors have converged respectively. The green dots show the effective number of SpMVs when we always use SpMM regardless of how many eigenvectors have converged.  As we can see from this figure,  the number of effective SpMVs is relatively high when we switch to RMM-DIIS after a few block Lanczos iterations. This is because it will take RMM-DIIS longer to converge if the initial eigenvector approximations produced from the block Lanczos iterations are not sufficiently accurate.  Regardless of when we switch from block Lanczos to RMM-DIIS, using SpMM throughout the RMM-DIIS algorithm appears to almost always yield the lowest effective SpMV count.    We can also see that the difference in the effective SpMV count is relatively large when we switch from block Lanczos to RMM-DIIS too early.  This is understandable because some of the approximate eigenvectors are more accurate than others when we terminate the block Lanczos iteration too early. As a result, the number of RMM-DIIS steps required to reach convergence may vary significantly from one eigenvector to another.  The difference in the rate of convergence prevents us from batching several SpMVs into a single SPMM. If we switch after more block Lanczos iterations have been performed, this difference becomes quite small as all approximate eigenvectors are sufficiently accurate and converge more or less at the same rate. The relative difference ($\tau$) between eigenvalue approximations from two consecutive RMM-DIIS iterations falls below $10^{-7}$ at the 13th block Lanczos iteration.  If we switch to RMM-DIIS at that point and use SpMM throughout the RMM-DIIS iteration, the number of effective SpMVs used in the hybrid scheme is about 40, which is slightly higher than the optimal 32 effective SpMVs required if we were to switch to RMM-DIIS after 23 block Lanczos iterations.

\begin{figure}[!tb]
\centering
\includegraphics[scale=0.4]{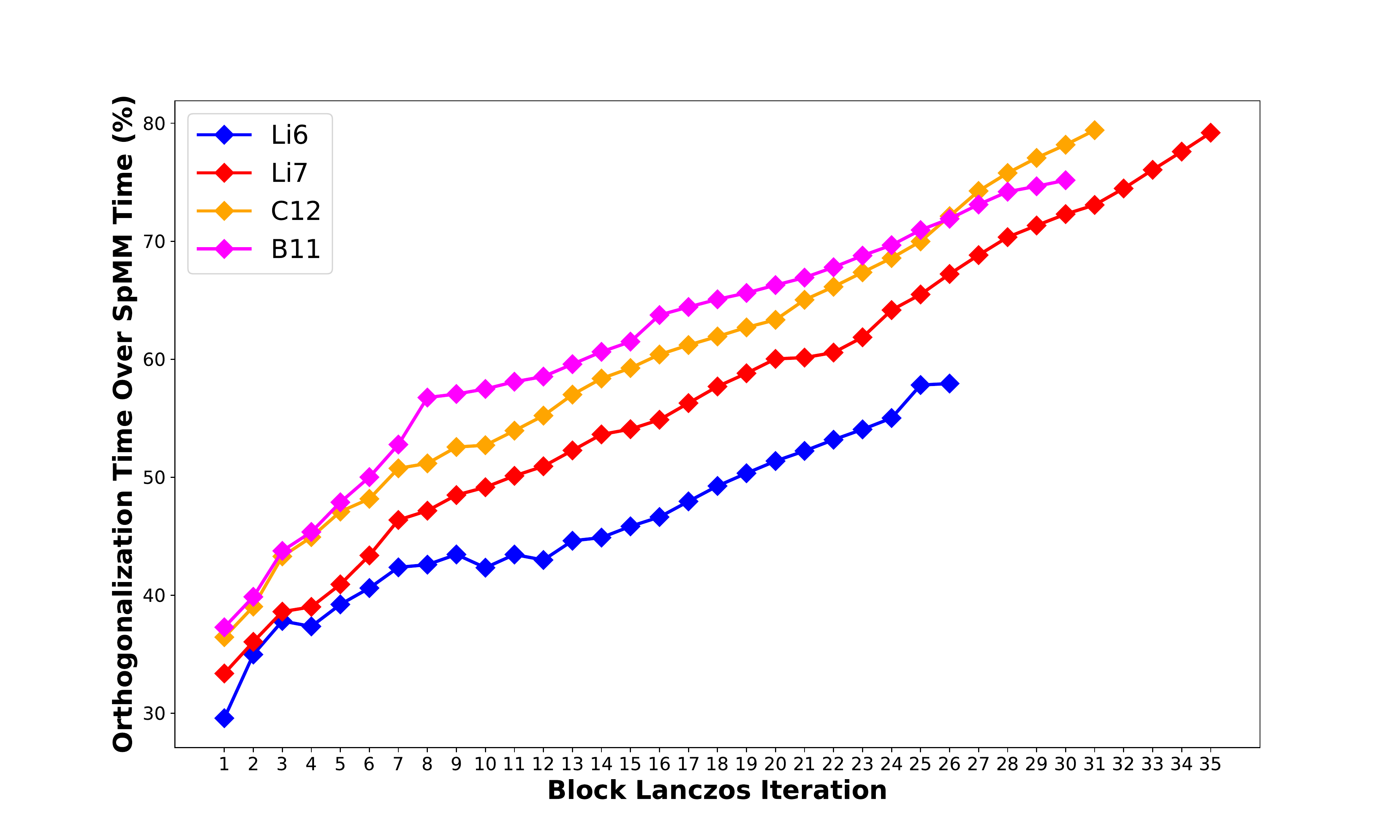}
\caption{The percentage of cumulative time spent on orthogonalization compared to SpMM in the block Lanczos algorithm for all test problems to compute the lowest five eigenvalues.}
\label{fig:BlockLanOrthCost}
\end{figure}
An early switch allows us to keep the basis orthogonalization cost of the block Lanczos algorithm as low as possible. As we can see from Figure~\ref{fig:BlockLanOrthCost}, the cost of orthogonalization as a percentage of the SpMM cost can become significantly higher as we perform more block Lanczos iterations. In particular, for all test problems, the orthogonalization cost exceeds 50\% of the SpMM cost after 20 block Lanczos iterations.  We note that the performance shown in Figure~\ref{fig:BlockLanOrthCost} is measured from the wallclock time of an implementation of the block Lanczos algorithm in the MFDn software executed on a single node of the Perlmutter using one MPI rank and 64 threads. %statements in this section, and not only to Fig. 9.  Also, Perlmutter needs to be better defined.  I presume you mean Perlmutter-CPU at NERSC, and it needs a reference.}

\begin{figure}[!tb]
\centering
\includegraphics[scale=0.6]{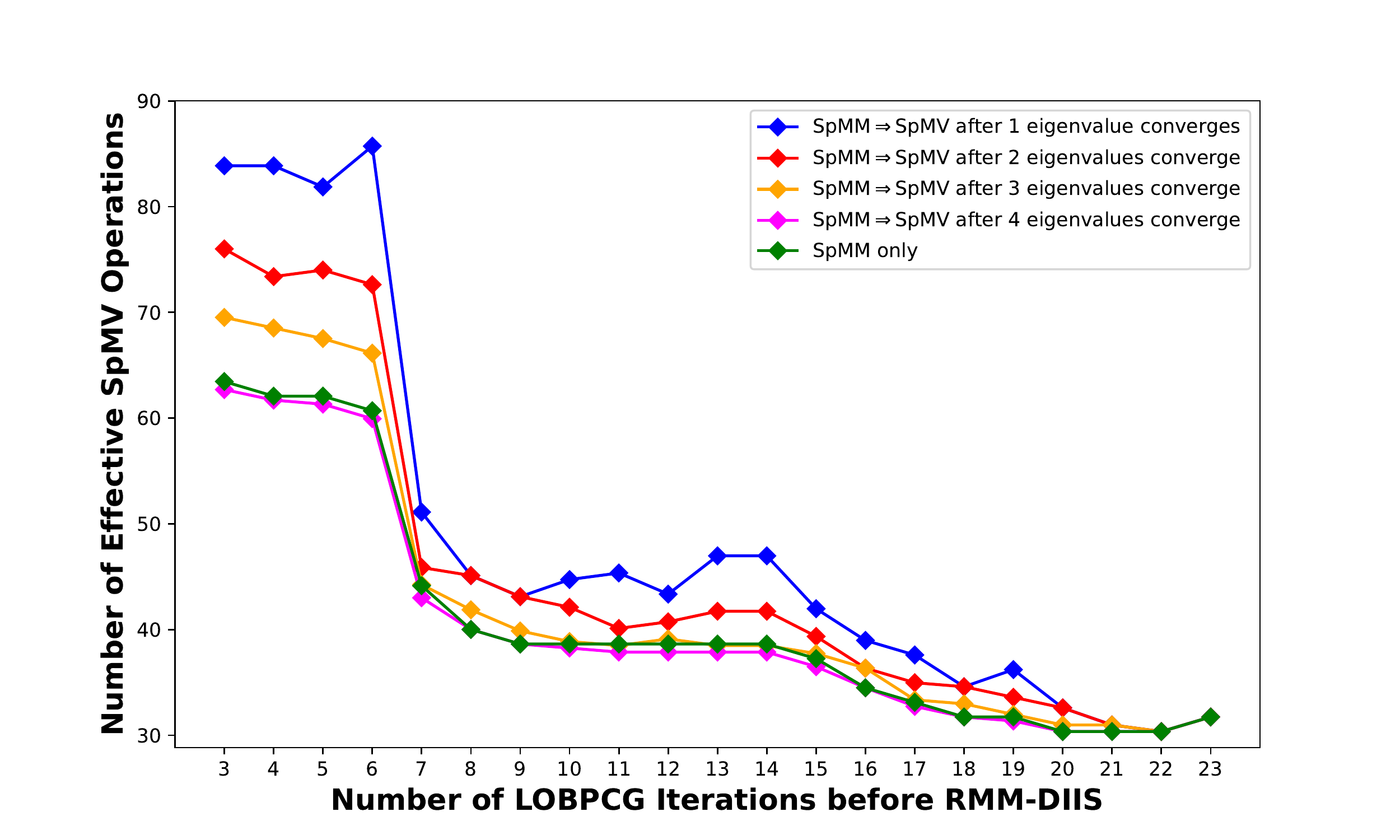}
\caption{Total effective SpMV cost of the LOBPCG/RMM-DIIS algorithm to compute the lowest five eigenvalues for $^6$Li w.r.t. switching strategy from SpMM to SpMV within RMM-DIIS.}
\label{fig:switch_lobpcg_li6}
\end{figure}
\begin{figure}[!tb]
\centering
\includegraphics[scale=0.4]{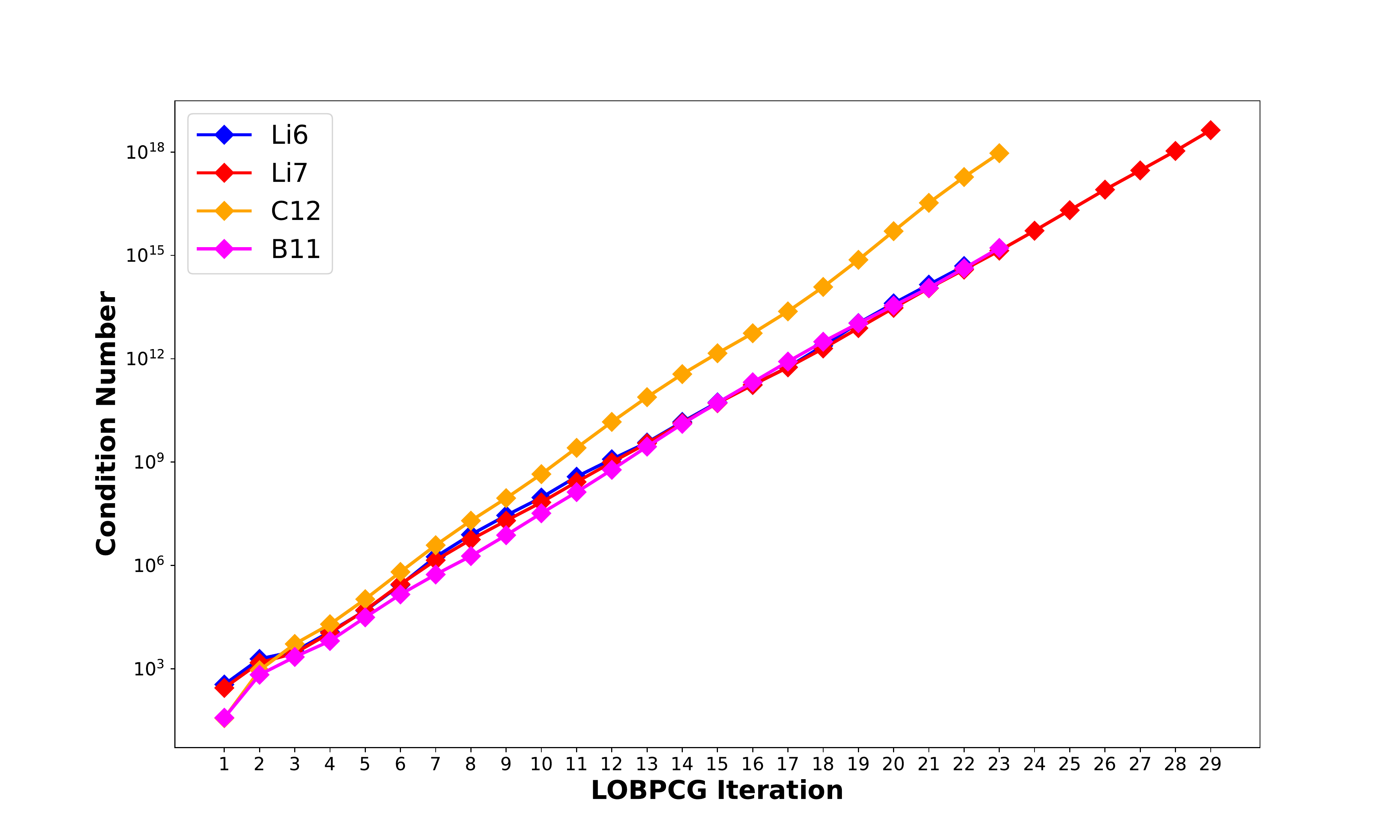}
\caption{The condition number of the LOBPCG subspace from which approximate eigenpairs are extracted to compute the lowest five eigenvalues of the $^6$Li Hamiltonian.}
\label{fig:Li6CondLOBPCG}
\end{figure}

Figure~\ref{fig:switch_lobpcg_li6} shows that it is beneficial to use SpMM throughout the hybrid LOBPCG/RMM-DIIS algorithm even after some of the eigenvectors have converged. When a sufficient number of LOBPCG iterations have been performed, the total SPMV count does not vary much.  In this case, we should terminate LOBPCG as early as possible to avoid potential numerical instabilities that can be introduced by the numerical rank deficiency of the preconditioned residual vectors~\cite{robustlobpcg}.   Figure~\ref{fig:Li6CondLOBPCG} shows that the estimated condition number of the subspace from which approximated eigenvalues and eigenvectors are extracted in the LOBPCG algorithm increases rapidly as we perform more LOBPCG iterations. Although the optimal SpMV count is attained when we switch to RMM-DIIS after 20 LOBPCG iterations, it is not unreasonable to terminate LOBPCG sooner, for example, after 12 iterations, when the estimated condition number of the LOBPCG subspace is around $10^{9}$. This is also the point at which the average relative change in the approximations to the desired eigenvalues $\tau$ just moves below $10^{-7}$. Therefore, the previously discussed strategy of using $\tau < 10^{-7}$ to decide when to switch to RMM-DIIS works well.  Even though this strategy would lead to a slight increase in the number of effective SpMV operations compared with the optimal effective SpMV count achieved when we switch to RMM-DIIS after 20 iterations, it makes the hybrid algorithm more robust and stable.  We should also note that in a practical calculation the optimal effective SpMV count and when the optimality is achieved is unknown a priori. This optimality is problem dependent, and is also architecture dependent.  
\section{Conclusion}
In this paper, we examine and compare a few iterative methods for solving large-scale eigenvalue problems arising from nuclear structure calculations. We observe that the block Lanczos and LOBPCG methods are generally more efficient than the standard Lanczos method and the stand-alone RMM-DIIS method in terms of the effective number of SpMVs. The Chebyshev filtering based subspace iteration method is not competitive with other methods even though it requires the least amount of memory and has been found to be very efficient for other applications. 
We show that by combining the block Lanczos or LOBPCG algorithm with the RMM-DIIS algorithm, we obtain a hybrid solver that can outperform existing solvers. The hybrid LOBPCG/RMM-DIIS method is generally more efficient than block Lanczos/RMM-DIIS when a good preconditioner is available. The use of RMM-DIIS in the block Lanczos/RMM-DIIS hybrid algorithm allows us to limit the orthogonalization cost in the block Lanczos iterations.  In the LOBPCG/RMM-DIIS hybrid algorithm, the use of RMM-DIIS allows us to avoid the numerical instability that may arise in LOBPCG when the residuals of the approximate eigenpairs become small.  We discuss the practical issue of how to decide when to switch from block Lanczos or LOBPCG to RMM-DIIS. A strategy based on monitoring the averaged relative changes in the desired approximate eigenvalues has been found to work well. 
Although the RMM-DIIS method is a single vector refinement scheme, we show the SpMVs in multiple independent RMM-DIIS iterations targeting different eigenpairs can be batched together and implemented as a single SpMM.  Such a batching scheme significantly improves the performance of the hybrid solver and is found to be useful even after some of the approximate eigenpairs have converged.

In large-scale nuclear structure calculations, we are primarily interested in the lowest eigenvalues and corresponding quantum numbers. In this case, we can use approximate convergence of the eigenvalues as a stop criterion, and write out the obtained eigenvectors.  Subsequently, one can then refine a subset of those eigenvectors (not necessarily those corresponding to the lowest eigenvalues) with the RMM-DIIS algorithm, before using those refined eigenvectors to evaluate additional observables of interest such as charge radii, magnetic and quadrupole moments, and electroweak transitions.

%\AB{In the last meeting we were talking about whether we want to explain what we anticipate might be happening at scale here.}
%\input appendix

\section*{Acknowledgements}
This material is based upon work supported by the U.S. Department of Energy, Office of Science, Office of Nuclear Physics under Grants DE-SC0023495 and DE-SC0023175, Office of Advanced Scientific Computing Research, Scientific Discovery through Advanced Computing (SciDAC) program via the FASTMath Institute under Contract No. DE-AC02-05CH11231, and National Science Foundation Office of Advanced Cyberinfrastructure under Grant 1845208. This work used resources of the National Energy Research Scientific Computing Center (NERSC) using NERSC Award ASCR-ERCAP m1027 for 2023, which is supported by the Office of Science of the U.S. Department of Energy under Contract No. DE-AC02-05CH11231.

\bibliographystyle{elsarticle-num}
\bibliography{refs_updated}

\end{document}